\documentclass[12pt,leqno]{amsart}
\usepackage{amssymb,amsmath}
\usepackage{latexsym, amssymb, amsmath, amsfonts}

in \textwidth=6.25 true in \headheight=.25 true
in \headsep=0.5 true cm \topmargin=0 true in
\begin{document}
\newcommand{\matr}[4]{\left[ \begin{array}{cc}
                              #1 & #2 \\
                              #3 & #4
                              \end{array} \right]}
\newcommand{\ra}{\rightarrow}
\newcommand{\Ra}{\Rightarrow}
\newcommand{\la}{\leftarrow}
\newcommand{\da}{\downarrow}
\newcommand{\lra}{\longrightarrow}
\newcommand{\lla}{\longleftarrow}
\newcommand{\lara}{\leftrightarrow}
\newcommand{\sra}[1]{\stackrel{#1}{\ra}}
\newcommand{\slra}[1]{\stackrel{#1}{\lra}}
\newcommand{\sla}[1]{\stackrel{#1}{\la}}
\newcommand{\slla}[1]{\stackrel{#1}{\lla}}
\newcommand{\sea}{\searrow}
\newcommand{\swa}{\swarrow}
\newcommand{\ov}{\overline}

\newcommand{\End}{{\rm End}\,}
\newcommand{\Gal}{{\rm Gal}\,}
\newcommand{\GL}{{\rm GL}\,}
\newcommand{\Gr}{{\rm Gr}\,}
\newcommand{\SL}{{\rm SL}\,}
\newcommand{\PSL}{{\rm PSL}\,}
\newcommand{\GSp}{{\rm GSp}\,}
\newcommand{\Aut}{{\rm Aut}\,}
\newcommand{\Ker}{{\rm Ker}\,}
\newcommand{\Prym}{{\rm Prym}\,}
\newcommand{\Jac}{{\rm Jac}\,}
\newcommand{\Pic}{{\rm Pic}\,}
\newcommand{\Sym}{{\rm Sym}\,}
\newcommand{\St}{{\rm St}}
\newcommand{\Id}{{\rm Id}}
\newcommand{\Tr}{{\rm Tr}}
\newcommand{\Disc}{{\rm Disc}\,}
\newcommand{\Mat}{{\rm Mat}\,}
\newcommand{\Nm}{\operatorname{Nm}}
\newcommand{\rank}{{\rm rank}\,}
\newcommand{\univ}{{\rm univ}}
\newcommand{\we}{{\rm w}}
\newcommand{\sing}{{\rm sing}}
\newcommand{\Per}{\operatorname{Per}}
\newcommand{\Br}{{\rm Br}}
\newcommand{\Shim}{{\rm Shim}}

\newcommand{\snu}{\mbox{$\scriptstyle{\nu}$}}
\newcommand{\snutil}
         {\mbox{$\scriptstyle{\til{\nu}}$}}
\newcommand{\spi}{\mbox{$\scriptstyle{\pi}$}}
\renewcommand{\sf}{\mbox{$\scriptstyle{f}$}}
\newcommand{\sg}{\mbox{$\scriptstyle{g}$}}
\newcommand{\sh}{\mbox{$\scriptstyle{h}$}}
\newcommand{\Imom}{\Im\omega}\renewcommand{\deg}{{\rm deg}\,}

\newcommand{\qi}{{\hat{\imath}}}
\newcommand{\qj}{{\hat{\jmath}}}
\newcommand{\qk}{{\hat{k}}}
\newcommand{\qmin}{\hat{\varepsilon}}

\newcommand{\quat}{G}
\newcommand{\til}{\tilde}
\newcommand{\Cpm}{{C_\pm}}
\newcommand{\al}{\alpha}
\newcommand{\be}{\beta}
\newcommand{\ga}{\gamma}
\newcommand{\lam}{\lambda}
\newcommand{\ip}[2]{\langle #1,\, #2 \rangle}
\newcommand{\iptil}[2]{{\langle #1,\,
                           #2 \rangle}_{\Ctil}}
\newcommand{\Ztwo}{\ZZ/2\ZZ}
\newcommand{\cTb}{\ov{{\cal T}}}
\newcommand{\taub}{\ov{\tau}}
\newcommand{\bmu}{{\bbB \mu}}

\newcommand{\bbb}[1]{{\bf #1}}
\newcommand{\bbB}{\mathbf}
\newcommand{\Aa}{{\bbB A}}
\newcommand{\CC}{{\bbB C}}
\newcommand{\FF}{{\bbB F}}
\newcommand{\PP}{{\bbB P}}
\newcommand{\QQ}{{\bbB Q}}
\newcommand{\RR}{{\bbB R}}
\newcommand{\ZZ}{{\bbB Z}}
\newcommand{\GG}{{\bbB G}}

\newcommand{\cA}{{\mathcal A}}
\newcommand{\cB}{{\mathcal B}}
\newcommand{\cC}{{\mathcal C}}
\newcommand{\cD}{{\mathcal D}}
\newcommand{\cH}{{\mathcal H}}
\newcommand{\cL}{{\mathcal L}}
\newcommand{\cM}{{\mathcal M}}
\newcommand{\cO}{{\mathcal O}}
\newcommand{\cP}{{\mathcal P}}
\newcommand{\cQ}{{\mathcal Q}}
\newcommand{\cR}{{\mathcal R}}
\newcommand{\cS}{{\mathcal S}}
\newcommand{\cT}{{\mathcal T}}
\newcommand{\cZ}{{\mathcal Z}}
\newcommand{\bM}{{\bf M}}
\newcommand{\bP}{{\bf P}}

\newcommand{\Ctil}{\til{C}}
\newcommand{\Ttil}{\til{T}}
\newcommand{\Ptil}{\til{P}}
\newcommand{\nutil}{\til{\nu}}

\newcommand{\Gmr}{\GG_m^r}

\newcommand{\pard}[2]{\partial{#1}/\partial{#2}}

\newcommand{\Pry}{P}
\newcommand{\cPry}{\cP}

\newcommand{\Pis}{\Pi^\ast}
\newcommand{\Pisp}{\Pi^{\prime,\ast}}
\newcommand{\rsm}{{\rm sm}}

\newcommand{\ssub}{\scriptstyle{\subset}}
\newcommand{\ssup}{\scriptstyle{\supset}}
\newcommand{\sslash}{\scriptstyle{/}}
\newcommand{\seq}{\scriptstyle{=}}
\newcommand{\smin}{\scriptstyle{-}}
\newcommand{\stimes}{\times}
\newcommand{\bto}{\ssub \seq \sslash \smin\smin }
\newcommand{\btt}{\stackrel{\subset}
                          {\ssub} \sslash \ssub }
\newcommand{\btth}
                {\ssub \ssub \sslash \smin\smin }
\newcommand{\btf}{\ssup\!\ssub \sslash  \stimes }
\newcommand{\tto}{\stackrel{\subset}{\ssub}\seq
                             \sslash\ssub \smin }
\newcommand{\ttt}
        {\ssub\ssub\seq \sslash \smin\smin\smin }
\newcommand{\tttp}{\ssup\!\ssub\seq
                      \sslash     \stimes\smin  }
\newcommand{\ttoX}{\ssub \smin\smin }
\newcommand{\tttX}{\ssub \ssub }
\newcommand{\Hdown}[1]{H_1(#1,\ZZ)}
\newcommand{\Hup}[1]{H^1({#1},{\ZZ})}
\newcommand{\piunC}{\pi_1(C,v)}
\newcommand{\piunCB}{\pi_1(C\setminus\Br,v)}
\newcommand{\hal}{\hat{\alpha}}
\newcommand{\hbe}{\hat{\beta}}
\newcommand{\hga}{\hat{\gamma}}
\newcommand{\hde}{\hat{\delta}}
\newcommand{\qu}{\hat{u}}
\newcommand{\hv}{\hat{v}}
\newcommand{\vtil}{\tilde{v}}
\newcommand{\altil}{ \tilde{\alpha}}
\newcommand{\betil}{\tilde{\beta}}
\newcommand{\gatil}{\tilde{\gamma}}
\newcommand{\detil}{\tilde{\delta}}

\newcommand{\zg}{\ZZ G}
\newcommand{\ZG}{\zg}
\newcommand{\zvq}{\ZZ V_4}
\newcommand{\ZVQ}{\zvq}
\newcommand{\chain}{\cC_\centerdot}
\newcommand{\chainmtag}{\chain(\cM')}
\newcommand{\chainzg}{\chain(\ZG)}
\newcommand{\chainvq}{\chain(\ZVQ)}
\newcommand{\cycle}{\cZ}
\newcommand{\ung}{1_G}
\newcommand{\minung}{(-1)_G}
\newcommand{\qib}{\ov{\imath}}
\newcommand{\qjb}{\ov{\jmath}}
\newcommand{\qkb}{\ov{k}}
\newcommand{\sigdiag}[4]
{\left(#1 \,\begin{array}{r} #2 \\ #3
\end{array}\, #4 \right)}

\newtheorem{lemma}{Lemma}
\newtheorem{proposition}[lemma]{Proposition}
\newtheorem{theorem}[lemma]{Theorem}
\newtheorem{claim}[lemma]{Claim}
\newtheorem{corollary}[lemma]{Corollary}
\newtheorem{remark}[lemma]{Remark}
\newtheorem{remarks}[lemma]{Remarks}
\newtheorem{definition}[lemma]{Definition}

\newenvironment{pf}{\noindent {\it Proof:}\/ }
               {\par \medskip \par}

\title[Abelian Varieties with Quaternion Multiplication]
{Abelian varieties with quaternion
multiplication}

\author{Ron Donagi}
\address{Mathematics department\\
University of Pennsylvania\\
Philadelphia, PA 19104-6395\\
USA\\}
\thanks{The first author thanks the Hebrew University
of Jerusalem for its hospitality. He was also
partially supported by NSF grants DMS-0104354
and FRG-0139799.}
 \email{donagi@math.upenn.edu}
\hfill
\author{Ron Livn\'e}
\address{Mathematics Institute\\
The Hebrew University of Jerusalem\\
Givat-Ram, Jerusalem  91904  Israel}
\email{rlivne@sunset.ma.huji.ac.il}
\thanks{The second author thanks the University
of Pennsylvania for its hospitality during fall
1997, when the main ideas of this work were
worked out}
\thanks{This work was partially supported by an
Israel-US BSF grant}
\date{June 9, 2005}
\keywords{Shimura varieties, quaternions, Prym
varieties} \subjclass{Primary: ; Secondary: }
\begin{abstract}

In this article we use a Prym construction to
study low dimensional abelian varieties with an
action of the quaternion group. In special cases
we describe the Shimura variety parameterizing
such abelian varieties, as well as a map to this
Shimura variety from a natural parameter space
of quaternionic abelian varieties. Our
description is based on the moduli of cubic
threefolds with nine nodes, a subject going back
to C. Segre, which we study in some detail.
\end{abstract}

\maketitle

\section{Introduction} In this article we study
the moduli of low dimensional abelian varieties
on which the quaternion group $G$ acts.  These
abelian varieties with quaternion multiplication
are realized as Prym varieties associated to
$G$\/-Galois covers $\Ctil/C$ of curves.  The
Prym construction gives a map $\Phi:\cM\ra\Shim$
from an appropriate parameter space $\cM$ to a
Shimura variety $\Shim$. Explicitly, $\cM$
parameterizes triples $(C,\Br,{\rm type})$
consisting of a curve $C$, a finite subset $\Br$
over which $\Ctil$ is ramified, and the type of
the ramification of $\Ctil$ above each point of
$\Br$. For dimension reasons, there are only
five cases in which the resulting map $\Phi$ can
be surjective. These are listed in
(\ref{cases}). The two cases of unramified
quaternion covers, which occur over a base curve
$C$ of genus $2$ or $3$, were considered
independently by van Geemen and Verra \cite{GV}.
Their work concentrated mainly on the genus $3$
case, while we study mainly the genus $2$ case.
The intersection of our results and theirs is
therefore quite small, and is indicated in
Remarks~\ref{GV1} and~\ref{GV2}. In the
unramified genus $2$ case and two ramified
cases, related to it through a degeneration, we
go further and completely determine the relevant
Shimura variety, which turns out to be the
modular curve $Y_0(2)/w_2$ (see
Corollary~\ref{PEL2}).

In a second part of our work we give a
relationship between our Prym varieties and
cubic threefolds with nine nodes, which extends
Varley's treatment (see \cite{var, don}) of the
ten nodal Segre cubic threefold. In particular,
our study of nine nodal cubics gives that their
moduli space is canonically the same modular
curve. In Theorem~\ref{moduli} we describe
several related moduli problems, a description
which might be of independent interest (see e.g.
\cite{CLSS}).  This allows us, in
Theorem~\ref{main}, to determine the fibers of
our map when $\Ctil$ is an unramified cover of a
curve $C$ of genus $2$. We also get the
surjectivity of $\Phi$ in our case and in the
two other related cases (Corollary~\ref{end}).
We thank I. Dolgachev for the reference to the
classical work of C. Segre (\cite{seg}) on cubic
threefolds with 9 nodes.

Our interest in the problem was raised by a
question of W. Baily, connected with his attempt
to find a moduli problem pertaining to the
Pl\"ucker embedding of the Grassmanian $G(2,6)$
into $\PP^{14}$ (the third exceptional domain
in~\cite{Bai}). In email correspondence from
1997--1998 he suggested the following problem.
Let $T^{\ast\ast} \ra T$ be a cyclic unramified
four sheeted cover of a genus $7$ trigonal curve
$T$, and identify $\Gal(T^{\ast\ast} / T)$ with
$\left<\qi\right>$ (here $\qi\in G$ is the
standard quaternion of order $4$\/). Let
$T^\ast/T$ be the intermediate $2$\/-sheeted
cover. The Prym variety
$\Prym(T^{\ast\ast}/T^\ast)$ is principally
polarized, it has dimension $12$, and it admits
an action by $\left<\qi\right>$. Baily remarked
that the space of such covers $T^{\ast\ast}/T$
has dimension $12$, as does the Shimura variety
parameterizing $12$\/-dimensional ppav's with a
$G$\/-action. Based on this and other
considerations, he asked whether
$\Prym(T^{\ast\ast}/T^\ast)$ admits an action by
$G$ (extending the action by
$\left<\qi\right>$). We sketch our (negative)
solution in the Appendix. It illustrates the
paucity of means for producing abelian varieties
with quaternion action by geometric means, and
thus indicates that the direct method in our
paper is probably difficult to circumvent. It
also proves a result which must be well-known
(compare \cite{Poo,KS}), namely that the
endomorphism ring of a generic hyperelliptic
jacobian is $\ZZ$.

\section{The basic construction}
\label{prel} The quaternion group $\quat$ has
order $8$ and a center of order $2$ generated by
$\qmin$. The standard generators $\qi$, $\qj$
for $\quat$ satisfy $\qi^2 = \qj^2 = \qmin$ and
$\qk :=\qi\qj = \qmin\qj\qi$. It has a unique
irreducible complex representation on which
$\qmin$ acts as $-1$. This representation $\St$
is $2$\/-dimensional. It is not defined over
$\RR$, but its sum with itself has a model $B$
over $\QQ$, which is unique up to an
isomorphism. In terms of the group ring $\QQ G$
we may take $B=\QQ G/(1+\qmin)\QQ G$. The
quotient of $\quat$ by $\left<\qmin\right>$ is
the Klein four-group $V_4 \simeq (\ZZ/2\ZZ)^2$.

We shall study $\quat$\/-Galois covers of
compact Riemann surfaces (``curves'') $\Ctil \ra
C$. Putting $\Cpm = \Ctil/ \left<\qmin\right>$,
and $C_t = \til{C}/\left<t\right>$ for $t =
\qi$, $\qj$, or $\qk$, we obtain a tower
\begin{equation}
\label{E:tower}
\begin{array}{rclcl}
&&\til{C}&&\\
&&\da\pi&& \\
&&\Cpm&&\\
&\swa &\da& \sea\\
C_\qi && C_\qk && C_\qj\\
&\sea & \da & \swa& \\
&&C\,.&&
\end{array}
\end{equation}
$\quat$ acts on the Prym variety $\Pry =
\Prym{\til{C}/\Cpm}$, with $\qmin$ acting as
$-1$. Hence $H_1(\Pry,\QQ)$ is a sum of copies
of $B$ and in particular $\Pry$ is even
dimensional.

As we shall see in Section~\ref{S:shim}, $\Pry$
has a natural PEL structure in the sense of
\cite{shi}. Hence we get a map
$\Phi:\cM\ra\Shim$ from the parameter space
$\cM$ of quaternion covers $\Ctil/C$, into an
appropriate moduli space of abelian varieties
with PEL structure, which is a Shimura variety
$\Shim$.
\begin{lemma} \label{types}
Let $g = g(C)$ be the genus of $C$, and suppose
there are $a$ branch points of $C$ over which
$\Ctil$ is ramified. If our map $\Phi$ is
surjective, with $\dim\Shim>0$, then $(g,a)$
must be $(0,4)$, $(1,2)$, $(1,3)$, $(2,0)$,
$(2,1)$ or $(3,0)$. \end{lemma} \begin{proof} If
the map is surjective, we must have $\dim\cM\geq
\dim\Shim >0$ by our assumption. Clearly there
are $3g-3+a$ parameters. On the other hand,
suppose that $C_\pm$ is unramified over exactly
$a'$ of the branch points --- we'll refer to
such branch points as being of the first type,
and to the others (over which $C_\pm$ is
ramified) as being of the second type. The
stabilizers in $\quat$ of ramification points
must be cyclic, necessarily $\pm 1$ over points
of the first type and of order $4$ over points
of the second type. By the Riemann-Hurwitz
formula we get
\[g(C_\pm) = 4g-3+a-a' \qquad \text{and} \qquad
  g(\Ctil) = 8g-7+3a-a',\]
and hence
\[\dim\Pry = g(\Ctil) - g(C_\pm) =
4(g-1) + 2a.\]

In loc.\ cit.\ Shimura shows that $\dim\Shim =
n(n-1)/2$ with $n = \dim\Pry/2$. Our assumption
$\dim\Shim > 0$ implies $2g-2+a=n\geq 2$.
Moreover the surjectivity of $\Phi$ implies the
second condition
\[3g-3+a \geq (2(g-1)+a)(2(g-1)+a-1)/2,\]
or equivalently $4\geq (2g+a-4)^2 + a = (n-2)^2
+ a$, so $0\leq a \leq 4$. We also have $n-2\geq
0$, and these two conditions have the solutions
indicated.
\end{proof}
Concerning the polarization, we have the
following
\begin{lemma}
The family of Pryms $\cPry$ is naturally
principally polarized if $a=0$. It is naturally
isogenous to a principally polarized family
$\cPry'$ if $a' = 0$.
\end{lemma}
\begin{proof}
The first part is well-known.
By~\cite[Lemma~1]{DL}, to get principally
polarized varieties which are isogenous to the
$\Pry$\/'s, we must make $C_\pm$ and $\Ctil$
singular by identifying the ramification points
in pairs. To do so $G$\/-equivariantly is
possible only if $a' = 0$, and only if we
identify the two points above each branch point
of $C$ of the second type.
\end{proof}

For $t=\qi$, $\qj$, or $\qk$, suppose there are
$a_t$ points of the second type above which
$C_t$ is unramified. We get a disjoint partition
of the branch locus $\Br$ corresponding to $a =
|\Br| = a' + a_\qi + a_\qj + a_\qk$.
\begin{lemma} \label{a_t}
The $a_t$\/'s all have the same parity. For each
genus $g$ and four non-negative integers $a'$,
$a_t$, with all the $a_t$\/'s of the same
parity, the quaternion towers with these
invariants form a complex space $\cM =
\cM(g;a',a_\qi,a_\qj,a_\qk)$ which is
quasi-projective.
\end{lemma}
\begin{proof}
Most of the claims follow directly from the
Riemann-Hurwitz formula. The parameter spaces
are quasi-projective since they are finite
covers of $\cM_g$.
\end{proof}
Lemma \ref{a_t} shows that exactly five of the
six cases allowed by Lemma \ref{types} actually
occur: we have found five principally polarized
cases with quaternion action for which
$\dim\cM\geq \dim \Shim>0$. In what follows we
restrict to these cases. Since  $a'=0$, we
change our notation from
$\cM(g;0,a_\qi,a_\qj,a_\qk)$ to
$\cM(g;a_\qi,a_\qj,a_\qk)$. In the following
table we give only the cases for which $a_\qi
\geq a_\qj \geq a_\qk$.
\begin{equation}
\label{cases}
\begin{array}{cccccc}
& g(C) &  \dim\Pry & \dim\cM & \geq & \dim\Shim\\
\cM(0,2,2,0) & 0  & 4 & 1& & 1\\
\cM(1,2,0,0) & 1  & 4 & 2& & 1\\
\cM(2,0,0,0) & 2  & 4 & 3& & 1\\
\cM(1,1,1,1) & 1  & 6 & 3& & 3\\
\cM(3,0,0,0) & 3  & 8 & 6& & 6
\end{array}
\end{equation}

In \cite{GV} van Geemen and Verra proved that
the spaces $\cM(g,0,0,0)$ were connected. We
will generalize this to the general case. We
have chosen to redo the unramified cases because
we need the general set-up as well as the
explicit forms of the corresponding
homomorphisms of fundamental groups to prove the
last part of the following Proposition:
\begin{proposition}
\label{strata} The five spaces
$\cM(g;a_\qi,a_\qj,a_\qk)$ above are connected
and irreducible. The space $\cM(1;2,0,0)$ can be
viewed as a boundary component of
$\cM(2,0,0,0)$, so its points represent curves
which are degenerations of  curves represented
by points of $\cM(2,0,0,0)$. Similarly,
$\cM(0;2,2,0)$ can be viewed as a boundary
component of $\cM(1;2,0,0)$.
\end{proposition}
\begin{proof}
For a given $C$ and $\Br$ (with its partition),
the tower \eqref{E:tower} with the $G$\/-action
is equivalent to a surjective homomorphism
$\psi: \piunCB \ra \quat$ up to an inner
automorphism. (Since we are working up to inner
automorphisms, the choice of the base point
$v\in C$ is unimportant.) This $\psi$ must send
each loop $\ga_i$ around a branch point
$P\in\Br$ to an element of order $4$ in $G$:
since $C_\pm/C$ is ramified at $P$, the
composition $\psi'$ of $\psi$ with the quotient
map to $V_4$ sends $\ga_i$ to an element $\neq
\ov{1}\in V_4$. Hence each of our $\cM$\/'s is a
finite cover of the {\em connected, irreducible}
moduli space of possible $(C,\Br)$. To prove
$\cM$ is connected and irreducible, it will
therefore suffice to show that for a fixed
$(C,\Br)$ the different $\psi'$\/'s are
contained in the image in $\cM$ of an
irreducible variety.

For this, let $\pi$ be the group defined by
generators $\al_i$, $\be_i$, for $i=1,\dots,g$,
and $\ga_j$, $j=1,\dots,a=|\Br|$, with one
defining relation
\[[\al_1,\be_1]\dots[\al_g,\be_g]
   \ga_1\dots\ga_a =1.\]
We next present $(C,\Br)$ in the usual way: we
view $C$ as a $4g$\/-sided polygon with sides
identified together in pairs. We moreover choose
non-intersecting paths from the base point $v$,
which we put at a corner of the polygon, to the
``punctures'', contained in the interior of the
polygon. It is well-known that this gives a
``standard'' isomorphism $\mu$ of $\pi$ with
$\piunCB$. The data $(C,\Br,\mu)$, where $\mu$
is such a standard isomorphism, given up to an
inner automorphism, is the same as a point in
the corresponding Teichm\"uller space $\cT =
\cT(g,|\Br|)$. This is the usual construction of
$\cT$, which is known to be irreducible (hence
connected) as a cover of $\cM$. (In itself this
is not enough to prove that $\cM$ is connected,
as several copies of $\cT$ may be needed to
cover it.) From each $\mu$ we get others
belonging to $\cT$ by elementary moves
consisting of Dehn twists and braiding. These
generate the mapping class group in $\Aut\cT$
(on which inner automorphisms act trivially).

In the sequel it will be convenient to identify
$\pi$ with $\piunCB$ via some standard
isomorphism. Then if $\psi,\psi':\pi\ra G$ as
above differ by an element $\nu$ of the mapping
class group (namely $\psi' = \psi\nu$\/), it
follows that the points they represent lie in
the same component of $\cM$. To prove that $\cM$
is connected and irreducible it thus suffices to
show that any two $\psi$\/'s differ by an
element of the mapping class group. We will
achieve this by a reduction to abelian subgroups
and quotients of $G$.

The point of reducing to the abelian case is to
be able to use the following well-known fact.
The action of the mapping class group induced on
the abelianization ($\simeq \ZZ^{2g+|\Br|}$\/)
of $\piunCB$, surjects onto the group which
permutes the $\gamma_j$\/'s and which sends each
$\al_i$ to $\sum_k n_{i,k}\al_k +n_{i,g+k}\be_k
+ \sum_l m_{i,l}\ga_l$ and each $\be_i$ to
$\sum_k n_{g+i,k}\al_k +n_{g+i,g+k}\be_k +
\sum_l m_{g+i,l}\ga_l$. Here the $2g\times 2g$\/
matrix $n_{i,j}$ is any integral matrix which is
symplectic for the cup product, and $n_{i,j}$,
$m_{i,j}$ are integers.

Using this fact we first bring to normal form
the composition $\psi_\pm:\pi\ra V_4$ of $\psi$
with the quotient map $G\ra V_4$. Since
$\psi_\pm$ factors through the canonical
quotient $H_1(C\setminus\Br,\ZZ/2\ZZ)$ of $\pi =
\piunCB$, we get by duality a copy $V_\pm$ of
$V_4$ in $H^1(C\setminus\Br,
\frac{1}{2}\ZZ/\ZZ)$.

Consider the unramified cases first. Then this
$H^1$ is the group $\Jac(C)[2]$ of points of
order $2$ of the jacobian
\[\Jac(C) = \Pic^0(C) = H^1(C,\cO_C)/H^1(C,\ZZ).\]
Algebro-geometrically, $V_\pm$ is the kernel of
the norm map $\Nm: \Jac(\Cpm) \ra \Jac(C)$. We
now have the following
\begin{lemma}
\label{liftpsi}
\begin{enumerate}
\item $V_\pm$ is isotropic for the Weil pairing
$w_2$ on $\Jac(C)[2]$.
\item Conversely, given a copy $V$ of $V_4$ in
$\Jac(C)[2]$, totally isotropic for the Weil
pairing, let $\psi_\pm: \piunC \ra V$ be the
corresponding homomorphism. Then $\psi_\pm$ can
be lifted to a homomorphism $\psi:\piunC \ra
\quat$. Two such lifts differ by multiplication
by an arbitrary homomorphism $\chi:\piunC \ra
\pm 1$; in particular, there are $16$ such
$\psi$\/'s.
\item For each lift $\psi$ let
$\Ctil_\psi\ra \Cpm$ be the corresponding cover.
Then there are $4$ inequivalent
$\Ctil_\psi$\/'s, and for each of them there are
$4$ actions of $G$ differing by inner
automorphisms.
\end{enumerate}
\end{lemma}
\begin{proof} {\em 1.}
The symplectic group over $\ZZ$ surjects onto
its mod $2$ analog, and the Weil pairing on
$\Jac(C)[2]$ ``is'' the mod $2$ cup product. By
Witt's theorem for symplectic forms we can bring
$V_\pm$, and hence $\psi_\pm$, to normal form by
a choice of some standard  $\mu$, so that we
have either

\noindent {\bf 1. }$V_\pm$ is isotropic for the
Weil pairing, and $\psi_\pm$ is given by
$\psi_\pm(\al_1) = \ov{\imath}$,
$\psi_\pm(\al_2) = \ov{\jmath}$,
$\psi_\pm(\al_i) = 1$ for $i\geq 3$, and
$\psi_\pm(\be_i) = 1$ for all $i$, with
$\ov{\imath}$, $\ov{\jmath}$ denoting the
respective images of $\qi$ and $\qj$ in $V_4$;
or

\noindent {\bf 2. } The Weil pairing is
non-degenerate on $V_\pm$, and $\psi_\pm$ is
given by $\psi_\pm(\al_1) = \ov{\imath}$,
$\psi_\pm(\be_2) = \ov{\jmath}$, and
$\psi_\pm(\al_i) = \psi_\pm(\be_i) = 1$ for all
$i \geq 2$.

However in case 2 we get that
$\psi([\al_i,\be_i])$ is $-1$ for $i = 1$ and is
$1$ otherwise, which is incompatible with the
defining relation of $\pi_1$. This shows that
case 1 must hold, proving the first part of the
lemma.

For the second part, suppose we are now in case
1 with $\psi_\pm$ in the above normal form. We
can then lift it to the {\em normal form}
$\psi_0:\piunC \ra \quat$ by setting
\begin{equation}
\label{E: standard}
 \psi_0(\al_1) = \qi,\;\;
\psi_0(\al_2) = \qj,\;\; \psi_0(\al_i) = 1\;\;
\mbox{\ for\ } i\geq 3 \mbox{\quad and \quad}
\psi_0(\be_i) = 1 \mbox{\ for all\ } i.
\end{equation}
That $\psi$ is well-defined and unique up to a
homomorphism to $\pm 1$ is clear. It is clear
that there are $16$ such homomorphism,
concluding the proof of the second part.

\noindent {\bf 3. } Finally, notice that the
image of $H^1(\Ctil,\ZZ/2\ZZ)\ra
H^1(\Cpm,\ZZ/2ZZ)$ is a copy of $V$, which has
order $4$. Hence there are exactly $4$ possible
$\Ctil$. Since conjugating a given $G$\/-action
on a given $\Ctil$ by an element of $G$ clearly
gives $4$ lifts of the same $V$\/-action of
$\Cpm$, we get the last part of the Lemma.

\end{proof}

We next show how to bring any lift $\psi$ to
normal form. We have $g=2,3$. Let $F_1$ be the
free group on generators $\al_1$, $\be_1$, and
let $F_2$ be the free group on $\al_i,\be_i$ for
$2\leq i\leq g$. These groups come with obvious
maps to $\pi$. (These maps are inclusions but we
will not use it except to omit their maps to
$\pi$ from the notation.) Notice that $F_1$ is
the fundamental group of a once punctured genus
$1$ surface, and that $F_2$ is the fundamental
group of a once punctured genus $g-1$ surface.
In addition $\psi(F_1)$ is cyclic of order $4$
on $\qi$ and $\psi(F_2)$ is cyclic of order $4$
on $\qj$. It is clear that each elementary
transformation $\tau$ on $F_1$ or $F_2$ ``is''
an elementary transformation on $\pi$. Viewing
the pair $(\psi(\al_1),\psi(\be_1))$ and the
$(2g-2)$\/-tuple
$(\psi(\al_2),\dots,\psi(\be_{2g}))$ as vectors
in $(\ZZ/4\ZZ)^2$, $(\ZZ/4\ZZ)^{2g-2}$
respectively, we may use a symplectic
transformation mod $4$ to move them to the first
unit vectors $e_1=(1,0,,\dots,0)$ of the
respective lengths $2$, $2g-2$. Lifting to the
mapping class group as before, we may get the
normal form $\psi(\al_1)=\qi$,
$\psi(\al_2)=\qj$, and the other generators map
to $0$. As was explained this implies that $\cM$
is irreducible and connected. (The proof did not
use the assumption $g=2,3$ so that $\cM$ is
irreducible and connected in general; however
the map to the Shimura variety is not surjective
for $g\geq 4$.)

We consider now the case of $\cM(1;1,1,1)$. To
put $\psi_\pm$ in normal form we apply it to the
defining relation of $\pi$ to get
$0+\psi_\pm(\ga_1) + \psi_\pm(\ga_2) +
\psi_\pm(\ga_3) = 0\in V_4.$

 As none of the $\psi_\pm(\ga_i)$\/'s may be $0$,
they must be a permutation of $\qib$, $\qjb$ and
$\qkb$. An appropriate element $h$ of the
mapping class group then allows us to permute
them so that $\psi_\pm(\ga_1)=\qib$,
$\psi_\pm(\ga_2)=\qjb$, and
$\psi_\pm(\ga_3)=\qkb$. In addition we may
assume that $\psi_\pm(\al_1)=\psi_\pm(\be_1)=0$
by taking an $h$ inducing an appropriate
translation modulo $2$. Now $\psi_\pm$ is in (a
unique) normal form. We next bring $\psi$ to (a
unique) normal form as well. First we use a
translation modulo $4$ to make
$\psi(\al_1)=\psi(\be_1)=0$, which is possible
as before since $\psi$ maps the group generated
by $\al_1$, $\be_1$ and $\ga_1$ to a cyclic
group of order $4$. We have $\psi(\ga_1)=\pm\qi$
and $\psi(\ga_2)=\pm\qj$, so using an inner
automorphism of $G$ we may assume both signs are
$1$, and then the defining relation forces
$\psi(\ga_3)=-\qk$. This proves the uniqueness
of a normal form and the irreducibility and
connectedness of $\cM(1;1,1,1)$ follows.

The remaining cases are similar. Straightforward
computations, whose details we omit, give the
normal forms $\psi(\al_1)=\qj$, $\psi(\be_1)=1$,
$\psi(\ga_1)=\psi(\ga_2)^{-1}=\qi$ for the
component $\cM(1;2,0,0)_\qi$ and
$\psi(\ga_1)=\psi(\ga_2)^{-1}=\qi$,
$\psi(\ga_3)=\psi(\ga_4)^{-1}=\qj$ for
$\cM(0;2,2,0)_\qk$. To get the other
possibilities we make a cyclic permutation on
$\qi$, $\qj$, and $\qk$. Notice that these are
(outer) automorphisms of $G$; the other outer
automorphisms --- those of order 2, such as
replacing $\psi(\al_1)=\qi$ by $\qk$ in the
normal form for $\cM(1;2,0,0)$, give equivalent
forms.

The statements regarding degeneration follow by
letting the curve acquire an ordinary double
point. Let $S(g,n)$ be a curve of genus $g$ with
$n$ punctures and with a base point $\ast$.
There is a standard inclusion $S(1,2)\subset
S(2,0)$ obtained by adding a handle connecting
the two punctures. In our standard presentations
for fundamental groups this corresponds to the
map $\pi_1(S(1,2),\ast) \ra \pi_1(S(2,0),\ast)$
given by $\ga_1\ra\al_2$ and $\ga_2 \ra
\be_2\al_2^{-1}\be_2$. Likewise, the map $\ga_3
\ra  \al_1$, $\ga_4 \ra \be_1\al_1^{-1}\be_1$
gives the map on fundamental groups
$\pi_1(S(0,4),\ast) \ra \pi_1(S(1,2),\ast)$
corresponding to the standard inclusion
$S(0,4)\subset S(1,2)$. Then our standard form
above for the map $\psi: \pi_1(S(2,0),\ast) \ra
G$ induces $\psi(\ga_1)=\qi$,
$\psi(\ga_2)=\qi^{-1}=-\qi$, $\psi(\ga_3)=\qj$,
and $\psi(\ga_4)=-\qj$, which are equivalent to
the normal forms for the two degenerate cases.
We omit the details. This completes the proof of
Proposition~\ref{strata}.
\end{proof}
\begin{remark}
\label{GV1} Van Geemen and Verra make a similar
construction for the special case that the cover
$\Ctil/C$ is unramified. They prove by a similar
method the connectedness of the parameter
spaces.
\end{remark}

\section{The Shimura varieties}
\label{S:shim} Let $B$ be the Hamilton
quaternion algebra over $\QQ$. It is generated
over $\QQ$ by $\quat$ with $\qmin$ going to
$-1$. Let $\Tr_{B/\QQ}: B\ra \QQ$ be the reduced
trace, defined by $\Tr(a_1+a_2\qi+a_3\qj+a_4\qk
= 2a_1$. Let $b\mapsto\bar{b}: B\ra B$ the main
involution (or conjugation) $\bar{b}=\Tr_{B/\QQ}
b-b$. The order $\bM'=\ZZ \langle
1,\qi,\qj,\qk\rangle$ is contained, with index
$2$, in a unique maximal order $\bM=\ZZ \qu
+\bM'$, where $\qu =(1+\qi+\qj+\qk)/2$. We have
$(\bM')^\times \simeq G$. (see \cite{vig}).

Set $P = \Prym(\til{C}/\Cpm)$ (see
diagram~\eqref{E:tower}). From the exponential
sheaf sequences $0\ra \ZZ \ra \cO \ra \cO^\times
\ra 0$ on $\Ctil$ and on $\Cpm$ we get an
identification
\begin{equation}
\label{E: identification}
 \Hdown{P} \simeq \left( \Ker \Hup{\Ctil} \sra{\Nm}
\Hup{\Cpm} \right) \simeq \left( \Ker
\Hdown{\Ctil} \sra{\pi_\ast} \Hdown{\Cpm}
\right)
\end{equation}
The polarization pairing $\ip{\,}{} : \Hdown{P}
\times \Hdown{P} \ra \ZZ $ is principal and
$2\ip{\,}{}$ is the restriction of the
intersection pairing on $\Hdown{\Ctil}$.

Recall that a polarization on an abelian variety
$A$ determines a Rosati involution $\rho$ on
$\End(A)\otimes \QQ$ (see \cite{mum}),
characterized by the property
\begin{equation}
\label{E:ros} \ip{mu}{v} = \ip{u}{\rho(m)v}
\end{equation}
for all $u$, $v$ in $H_1,A,\QQ)$ and $m\in
\End(A)\otimes\QQ$.

We have the following
\begin{lemma} \label{L: P}
The action of $G$ on $\Ctil$ induces an action
of $\bM'$ on $P$. The Rosati involution
preserves $\bM' \subset \End P$ and induces the
main involution on $\bM'$.
\end{lemma}
\begin{proof}
 From the definition of $P$ we get a $G$
action on $P$, with $\qmin$ acting as $-\Id_P$.
By linearity $\bM' \simeq \ZG/(1+\qmin)\ZG$ acts
on $P$. Since $G$ preserves the orientation on
$\Ctil$ it preserves the intersection pairing.
This implies the relationship $\ip{mu}{v} =
\ip{u}{m^{-1}v}$ for all $u,v\in \Hdown{P}$ and
$m \in \quat$. Since $m^{-1} = \ov{m}$ for any
$m\in \quat$, we get \eqref{E:ros} by linearity
for all $m\in \bM'$.
\end{proof}
We now specialize to the case when $g = 2$. Then
we have the following:
\begin{theorem}
\label{PEL} When $g=2$ and the
tower~\eqref{E:tower} is unramified we have the
following:
\begin{enumerate}
\item The $\bM'$ action on $P$ extends (uniquely)
to an $\bM$ action.
\item Under this action $\Hdown{P}$ is free of
rank $2$ over $\bM$, and
\item with respect to an appropriate $\bM$ basis
$\lambda_1$, $\lambda_2$ of $\Hdown{P}$ the
polarization pairing is given by
\begin{equation}
\label{E: pair} \ip{\sum_i m_i \lambda_i}{\sum_j
n_j\lambda_j}=\Tr_{B/\QQ} \sum_i \sum_j v_{ij}
\ov{m}_i n_j,
\end{equation}
 where $[v_{ij}]= \frac{1}{2}
\matr{2(\qi+\qj)}{-1-\qi}{1-\qi}{0}$.
\end{enumerate}
\end{theorem}
\begin{proof}
{\bf 1.\ } We must show that $\bM \Hdown{P}$, a
priori contained in $H_1(P,\QQ)\subset
H^1(\Ctil,\QQ)$, is in fact contained in
$\Hup{\Ctil}$. Since the space of towers
$\Cpm/C$ is connected, it suffices to verify the
inclusion $\bM H_1 \subset H_1$ for one $\Cpm/C$
(and all compatible $\Ctil$\/'s). Let $C$ be the
smooth projective model of $y^2=x^6-1$ and let
$P_\zeta$ denote the Weierstrass point
$(x,y)=(\zeta,0)$ of $C$, where $\zeta^6=1$. Let
$V\subset \Jac C[2]$ be the subgroup generated
by divisor classes of $P_\zeta + P_{-\zeta}
-K_C$ for $\zeta^3=1$. Then $\bmu_6$ acts on $C$
via $\zeta:(x,y) \mapsto (\zeta x,y)$, preserves
$V$, and induces on $V$ a $\bmu_3$-action. It
follows that $G$ and $\bmu_6$ generate in $\Aut
\Ctil$ a central extension $\tilde{A}_4$ of
$A_4=\bmu_3 \ltimes V$ by $\{\pm 1\}$. This
group is known to be the group of units of $\bM$
(\cite{vig}), and it spans $\bM$ additively.
Similarly ${\bM'}^\times=G$ spans $\bM'$
additively. Hence $\bM$ acts on $\Hdown{\Ctil}$
extending the $\bM'$ action.

\noindent {\bf 2.\ } Since $\bM$ has class
number $1$ (\cite{vig}) and $\Hup{P}$ is torsion
free, it follows that it is free, necessarily of
rank $2$ since $\dim P =4$.

\noindent {\bf 3.\ } As was already remarked,
$H_1(P,\QQ)$ is $B$\/-free since $-1\in \quat$
acts on $P$ as $-1$. As in the previous Lemma,
the $\quat$\/-action shows that the polarization
is $B$\/-skew-hermitian on $H_1(P,\QQ)$. Since
the trace form $(x,y)\in B^2 \mapsto
\Tr_{B/\QQ}\ov{x}y$ is non-degenerate, there
exist unique elements $v_{ij}\in B$ for
which~\eqref{E: pair} holds. The skew-symmetry
of the polarization implies that $v_{ji} =
-\ov{v_{ji}}$. (This part holds in general, not
just for the genus $2$ case, except that the
rank of $H_1(P,\QQ)$ over $B$ is usually not
$2$).

We will compute the pairing by an explicit (and
lengthy) calculation, in the course of which we
will in fact reprove parts (1) and (2).

Let $C$ have genus $2$ and let $\Ctil \ra C$ be
an unramified $\quat$-cover. Write
\[\piunC= \langle
\alpha,\beta,\gamma,\delta \mid
[\alpha,\beta][\gamma,\delta]=1 \rangle\] with
$v$ a base-point. As was shown we may (and do)
let $\phi: \piunC \ra \quat$ be the map
characterized by $\phi(\alpha)=\qi$,
$\phi(\gamma) =\qj$,
$\phi(\beta)=\phi(\delta)=1$. In the universal
cover $C^{\rm univ}$ of $C$ choose a base-point
$\hv$ above $v$ and lift $\alpha$, $\beta$,
$\gamma$, and $\delta$ to (not necessarily
closed) paths $\hal$, $\hbe$, $\hga$, $\hde$,
starting at $\hv$. Then $C^{\rm univ}$ is a copy
of $\RR^2$ subdivided into ``octagons'' by the
paths $\{\tau.\xi \mid \tau \in \piunC, \xi \in
\{ \hal, \hbe, \hga, \hde \}\}$.
\begin{figure}[!h] \caption{} \label{F: fund}
\setlength{\unitlength}{2mm}
\begin{picture}(50,23)(0,0)
\multiput(7,2)(0,17){2}{\line(7,0){7}}
\multiput(14,2)(-12,12){2}{\line(1,1){5}}
\multiput(19,7)(-17,0){2}{\line(0,1){7}}
\multiput(2,7)(12,12){2}{\line(1,-1){5}}
\put(7,2){\vector(1,0){3.5}}
\put(14,2){\vector(1,1){2.5}}
\put(19,14){\vector(0,-1){3.5}}
\put(14,19){\vector(1,-1){2.5}}
\put(14,19){\vector(-1,0){3.5}}
\put(7,19){\vector(-1,-1){2.5}}
\put(2,7){\vector(0,1){3.5}}
\put(7,2){\vector(-1,1){2.5}} \put(10.5,19.5)
{\makebox(0,0)[b]{$[\delta,\gamma]\hga$}}
\put(17,17)
{\makebox(0,0)[bl]{$[\alpha,\beta]\hbe$}}
\put(19.5,10.5)
{\makebox(0,0)[l]{$\alpha\beta\alpha^{-1}
\hal$}} \put(17,4){\makebox(0,0)[lt]{$\alpha
\hbe$}} \put(10.5,1.5){\makebox(0,0)[t]{$\hal$}}
\put(4,4){\makebox(0,0)[rt]{$\hde$}}
\put(1.5,10.5){\makebox(0,0)[r]{$\delta \hga$}}
\put(4,17){\makebox(0,0)[rb]
{$\delta\gamma\delta^{-1} \hde$}}
\put(7,1.5){\makebox(0,0)[t]{$\hv$}}
\put(7,2){\makebox(0,0){$\bullet$}}
\put(10.5,10.5){\makebox(0,0){$F$}}
\put(40,11){\makebox(0,0){$\bullet$}}
\put(33,11){\line(1,0){14}}
\put(35,6){\line(1,1){10}}
\put(40,4){\line(0,1){14}}
\put(45,6){\line(-1,1){10}}
\put(40,11){\vector(1,0){3.5}}
\put(40,11){\vector(1,1){2.5}}
\put(40,18){\vector(0,-1){3.5}}
\put(35,16){\vector(1,-1){2.5}}
\put(40,11){\vector(-1,0){3.5}}
\put(40,11){\vector(-1,-1){2.5}}
\put(40,4){\vector(0,1){3.5}}
\put(45,6){\vector(-1,1){2.5}}
\put(41.4,11.8){\makebox(0,0)[l]{$\hv$}}
\put(47.5,11){\makebox(0,0)[l]{$\hal$}}
\put(45.5,16.5){\makebox(0,0)[lb]{$\hde$}}
\put(40,18.5){\makebox(0,0)[l]{$\gamma^{-1}
\hga$}}
\put(34.5,16.5){\makebox(0,0)[br]{$\delta^{-1}
\hde$}} \put(32.5,11){\makebox(0,0)[r]{$\hga$}}
\put(34.5,5.5){\makebox(0,0)[tr]{$\hbe$}}
\put(40,3.5){\makebox(0,0)[t]{$\alpha^{-1}\hal$}}
\put(45.5,5.5)
{\makebox(0,0)[lt]{$\beta^{-1}\hbe$}}
\put(45,13){\makebox(0,0){$F$}}
\end{picture}
\end{figure}

The left part of  Figure~\ref{F: fund}  gives
the sides of the (unique) octagon $F\subset
\RR^2$ whose boundary contains both $\hal$ and
$\hde$. The right part gives a {\em planar}
neighborhood of $ \hv$.

Let $R$ be the group ring $\ZZ\piunC$, and let
$\chain$ be the chain complex
\begin{equation}\label{E:chain}
  0\ra C_2 \ra C_1\ra C_0 \ra 0
\end{equation}
where $C_i$ is the left $R$\/-free module on the
basis set $\{F\}$, $\{\hal,\hbe,\hga,\hde\}$,
and $\{\hv\}$ for $i=2$, $1$, and $0$
respectively.

The differentials $\partial_1$, $\partial_2$ are
defined using the right and left sides of
figure~\ref{F: fund}  respectively as follows:
\begin{align}
\partial_2(rF) &= r\left[
(1-\alpha \beta\alpha^{-1}) \hal +(\alpha
-[\alpha,\beta])\hbe +
([\delta,\gamma]-\delta)\hga + (\delta \gamma
\delta^{-1}-1)\hde \right]
\label{E: deltaF} \\
\partial_1(r\hat{\tau}) &= r(\tau-1) \hv
\label{E: deltatau}
\end{align}
for any $r\in R$, with $\tau$ denoting either of
the symbols $\alpha$, $\beta$, $\gamma$ or
$\delta$.

For any right $R$\/-module $M$ let $\chain(M)$
be the complex $M \otimes_R \chain$. For a basis
element $\sigma = F$, \dots, $\hv$ as above, we
will denote $1\otimes_R\sigma\in \chain(M)$
respectively by $F_M$, $\hal_M$, \dots,
$\hde_M$, $\hv_M$. If $M$ is an $R$\/-algebra,
these elements form a free $M$\/-basis for
$\chain(M)$, and we will write $mF_M$, \dots,
for $m\otimes_R F$, \dots respectively. We will
write $\cycle_i(M)$ for the $i$\/-cycles of
$\chain(M)$. The homology class of a cycle $\xi$
will be denoted by $[\xi]$ or simply by $\xi$.

In the special case $M=R$ we get back
$\chain(R)=\chain$, the cellular (or CW) chain
complex for $C^{\rm \univ}$. Since this
description is clearly $\piunC$\/-equivariant,
we get the chain complex for $\Gamma\backslash
C^{\rm \univ}$, for any subgroup $\Gamma \subset
\piunC$, by taking $M=\ZZ\left(\Gamma \backslash
\piunC\right)$. If $\Gamma$ is normal in
$\piunC$ then this description is
$\piunC/\Gamma$\/-equivariant. If now
$\Gamma=\Ker \phi$ then $\piunC/\Gamma \simeq G$
and the formulas \eqref{E: deltaF}, \eqref{E:
deltatau}, with $M=\ZG$ simplify to
\begin{align}
\partial_2(F_{\ZG})&= (\qi-1)\hbe_{\ZG} +
(\qj-1)\hde_{\ZG}
\label{E:deltaFsimp} \\
\label{E: deltatausimp}
    \begin{split}
    \partial_1(\hal_{\ZG})&= (\qi-1) v_{\ZG},
    \quad
    \partial_1(\hbe_{\ZG}) =
    \partial_1(\hde_{\ZG})= 0,
    \quad \text{and} \quad
    \partial_1(\hga_{\ZG})=(\qj-1) \hv_{\ZG}.
    \end{split}
\end{align}

Similarly, \eqref{E:chain} with $M=\ZVQ$ and
\eqref{E:deltaFsimp}, \eqref{E: deltatausimp}
also describe the chain complex for $\Cpm$ if we
replace $F_{\ZG}$, \dots, $\hv_{\ZG}$ by
$F_{\ZVQ}$, \dots, $\hv_{\ZVQ}$, and $\qi$,
\dots by their images in $\ZVQ$.

The order $\bM' \simeq \zg/(1+\qmin)\zg$ is
identified with $\Ker (\ZG \ra \ZVQ) = (1-
\qmin) \ZG$, by sending the image of $m\in \ZG$
in $\bM'$ to $(1-\qmin)m$. Since $\chain$ is
free, we see that the projection $\pi:\Ctil \ra
\Cpm$ induces a $G$-equivariant exact sequence
\begin{equation}
\label{E: kernel}
 0\ra \chainmtag \ra \chainzg \sra{\pi_*} \chainvq
\ra 0.
\end{equation}
We can (and will) therefore identify
$\chainmtag$ with $\Ker (\chainzg \sra{\pi_*}
\chainvq) = (1- \qmin) \chainzg $. As above, the
identification is between the image of $m\in
\chainzg$ in $\chainvq$ and $(1-\qmin)m\in
\chainmtag$. In particular, for the basis
elements we have $F_{\bM'} = (1-\qmin)F_{\ZG}$,
\dots, $\hv_{\bM'} = (1-\qmin)\hv_{\ZG}$. We now
have the following
\begin{proposition}
\label{P: L}
\begin{enumerate}
\item The torsion subgroup of $H_1(\chainmtag)$
has order $2$ and the quotient
$L=H_1(\chainmtag/{\rm torsion}$ is naturally
$\Hdown{P}$.
\item The natural $\bM'$ action on $L$ extends
(uniquely) to an $\bM$ action, and $L$ is
$\bM$\/-free on the classes of the cycles
$\lambda_1 = (\qi+1)\hal_{\bM'} - (\qj +1)
\hga_{\bM'}$ and $\lambda_2 = \hbe_{\bM'}$.
\end{enumerate}
\end{proposition}
\begin{proof}
\begin{enumerate}
\item The homology sequence of \ref{E: kernel}
gives the exact sequence
\[ 0\ra \ZZ/2\ZZ \ra H_1(\chainmtag) \ra
\Hdown{\Ctil} \ra \Hdown{\Cpm}, \]
 since $H_2(\Ctil,\ZZ)\stackrel{\pi_*}{\ra}
H_2(\Cpm,\ZZ)$ is identified with
$\ZZ\stackrel{\deg \pi}{\ra}\ZZ$, and $\deg\pi =
2$. As $\Hdown{\Ctil}$ is torsion-free, the
result follows from \eqref{E: identification}.
\item To analyze $H_1(\chainmtag)$ write
$C_1(\bM')=C' \oplus C''$, where $C'={\rm span}
\{\hal_{\bM'}, \hga_{\bM'}\}$ and $C''={\rm
span} \{\hbe_{\bM'},\hde_{\bM'}\}$. By
\eqref{E:deltaFsimp} and \eqref{E: deltatausimp}
we have
\[H_1(\chainmtag) = \Ker(\partial_1\mid C')
\oplus (C''/\partial_2(\bM'F)) \]
 We will show that
\begin{enumerate}
\item $\Ker(\partial_1|C')$ is preserved by
multiplication by $\bM$ and is $\bM$\/-free on
$\lambda_1$.
\item $\partial_2(\bM F)\subset C'' $ and the
$\bM'$\/-module structure on $C''/\partial_2(\bM
F)$ extends uniquely to an $\bM$\/-module
structure, for which $C''$\/ is $\bM$\/-free on
$\lambda_2$.
\end{enumerate}
\end{enumerate}
\end{proof}
These two assertions imply Proposition~\ref{P:
L}.2. In turn they follow from the following:
\begin{lemma} \label{L: A} Let $\chi: B^2 \ra B$
and $\psi: B\ra B^2$ be the maps of left
$B$-modules given by
\[ \chi(x,y)=x(\qi-1)+ y(\qj-1) \quad\mbox{ and }
\quad \psi(x)=x(\qi - 1, \qj-1)
 .\] Then
\begin{enumerate}
\item $(\bM' \oplus \bM') \cap \Ker \chi=\bM(\qi +1,
-\qj -1)$.
\item We have $(\bM' \oplus \bM') \cap \psi(B) =
\psi(\bM)$ and $(\bM' \oplus \bM')+ \psi(B) =
\bM(1,0) +\psi(B)$; the last sum is direct.
\end{enumerate}
\end{lemma}
\begin{proof} Observe the isomorphism of (commutative)
rings
\[\bM'/2\bM'\simeq \FF_2[\varepsilon,\varepsilon']
/\varepsilon^2=(\varepsilon')^2=0 \] given by
$1+\qi \mapsto \varepsilon$, $1+\qj \mapsto
\varepsilon'$. Then $2\qu = 1+\qi+\qj+\qk \in
\bM'$ maps to $1+ (1+\varepsilon)
           (1+\varepsilon')+(1+\varepsilon)
(1+\varepsilon')=\varepsilon \varepsilon'$. We
now prove the two assertions of the lemma.
\begin{enumerate}
\item Set $z=(\qi+1,-\qj-1)\in B^2$. Then
$\chi(z)=0$, so $Bz \subset \Ker \chi$. In
addition $uz$ maps to $\varepsilon
\varepsilon'(\varepsilon,\varepsilon')=(0,0)$ in
$(\bM'/2\bM')^2$, so $(u/2)\cdot z$ is in $\bM'
\oplus \bM'$. Hence
$\mbox{RHS}\subset\mbox{LHS}$. Conversely,
suppose $\chi(x,y)=0$. Then for
$t=x(\qi-1)=y(1-\qj)$ we have
$t(\qi+1,-\qj-1)=-2(x,y)$, so that $\Ker \chi
\subset Bz$. If in addition $(x,y) \in \bM'
\oplus \bM'$ then $t\in \bM'$ has reduction
$\ov{t}$ to $\bM'/2\bM'$ which is divisible both
by $\varepsilon$ and by $\varepsilon'$, hence
$\ov{t}$ is a multiple of
$\varepsilon\varepsilon'$ which is the reduction
of $2\qu$. It follows that $t\in 2\qu\ZZ +
2\bM'$, so that $t/2 \in \qu\ZZ + \bM' = \bM$ as
asserted.
\item Notice that the sum on the RHS is direct.
Since
\[\psi(1+j) \equiv(2\qu,0) \pmod {2\bM'\oplus 2\bM'} \]
we see that $(\qu,0)\in \mbox{ LHS}$, so that
$\mbox{RHS}\subset\mbox{LHS}$. Conversely,
$(\qj-1)^{-1} (\qi -1) \in \bM$, so that
 $(0,1) = \psi( (\qj-1)^{-1} ) -
   ((\qj-1)^{-1} (\qi-1),0)$ belongs to the RHS,
giving $\mbox{LHS}\subset\mbox{RHS}$. This
completes the proof of the Lemma and hence of
the Proposition.
\end{enumerate}
\end{proof}

We shall now compute the polarization form of
$P$ in terms of the $\bM$\/-basis $\lambda_1$,
$\lambda_2$ of $L$. For this we shall use the
embedding of $L$ into $H_1(\chainzg)$.
By~\eqref{E: pair} we have that
\[ \ip{\sum_i m_i \lambda_i}{\sum_j n_j
\lambda_j} = \sum_i \sum_j \Tr v_{ij}
\ov{m_i}{n_j}
\]
for any $m_i, n_i \in B$ and appropriate
$v_{ij}$'s. To determine them, we need only
determine $\ip{m\lambda_i}{\lambda_j}$ for any
$1 \leq i \leq j\leq 2$, and $m\in \{
1,\qi,\qj,\qk\}$. Denote the intersection
pairing on $H_1(\chainzg)\simeq \Hdown{\Ctil}$
by $\iptil{\, }{ }$. Then for any $x$,
$y\in\cycle_1(\ZG)$ mapping to $\ov{x}$, $\ov{y}
\in \cycle_1(\bM')$ we have the basic formula
\begin{equation}
\label{E: pairing}
\ip{\ov x}{\ov y}= \frac{1}{2}
\iptil{(1-\qmin)x}{(1 - \qmin)y} = \frac{1}{2}
\iptil{x}{(1 - \qmin)^2y} = \iptil{x}{(1 -
\qmin)y}.
\end{equation}
Indeed, the first equality is the definition of
the polarization pairing on $P$, the second
follows from the general formula
$\ip{gx}{gy}=\ip{x}{y}$ for all $g\in G$, and
the third holds since
 $(1 -\qmin)^2=2(1 - \qmin)$.

We will let $\vtil$, $\altil$, $\betil$,
$\gatil$ and $\detil$ denote the images of
$\hv$, $\hal$, $\hbe$, $\hga$ and $\hde$ in
$\Ctil$. \vspace{0.3cm}

\emph{Computation of $v_{22}$}\/:\ \ We are
assuming that $\phi:\piunC \ra G$ is in normal
form~\eqref{E: standard}. Thus $\betil$ is a
closed loop. We have
$\iptil{\hbe_{\ZG}}{\hbe_{\ZG}}=0$, and since
$g\betil$ and $\betil$ are clearly disjoint for
$g\in G-{1}$, we get
$\iptil{\hbe_{\ZG}}{g\hbe_{\ZG}}=0$ for all
$g\in G$. Therefore $\ip{\lambda_2}{\lambda_2} =
\Tr_{B/\QQ}v_{22}m =\iptil{\hbe_{\ZG}}{(1 -
\qmin)m\hbe_{\ZG}} = 0$ for all $m\in \bM'$,
whence $v_{22}=0$. \vspace{0.3cm}

\emph{Computation of $v_{11}$}\/:\ \ Lift
$\lambda_1$ to the cycle $\sigma = (1 + \qi)
\hal_{\ZG} - (1 + \qj)\hga_{\ZG} \in
\chain(\ZG)$, represented by the following:
\begin{figure}[!ht] \caption{} \label{F: loz}
\setlength{\unitlength}{2.5mm}
\begin{picture}(20,13)(0,0)
\put(4,7){\line(2,-1){6}}
 \put(4,7){\vector(2,-1){3}}
\put(4,7){\line(2,1){6}}
 \put(4,7){\vector(2,1){3}}
\put(10,10){\line(2,-1){6}}
 \put(10,10){\vector(2,-1){3}}
\put(10,4){\line(2,1){6}}
 \put(10,4){\vector(2,1){3}}
\put(3.5,7){\makebox(0,0)[r]{$\vtil$}}
\put(10,10.5){\makebox(0,0)[b]{$\qi \vtil$}}
\put(16.5,7){\makebox(0,0)[l]{$\qmin \vtil\;,$}}
\put(10,3.5){\makebox(0,0)[t]{$\qj \vtil$}}
\put(7,9){\makebox(0,0)[br]{$\altil$}}
\put(13,9){\makebox(0,0)[bl]{$\qi \altil$}}
\put(13,5){\makebox(0,0)[tl]{$\qj \gatil$}}
\put(7,5){\makebox(0,0)[tr]{$\gatil$}}
\end{picture}
\end{figure}

By \eqref{E: pairing} we need to compute
$\iptil{\sigma}{(1-\qmin)l\sigma}$ for $l=1$,
$\qi$, $\qj$, and $\qk$. Now for any $x \in
\cycle_1(\zg)$ and $g\in G$, $g\neq \pm 1$ we
have $\iptil{x}{\qmin x}= \iptil{\qmin
x}{\qmin{}^2 x} =-\iptil{x}{\qmin x}$. Hence
$\iptil{x}{\qmin x}=0$. Since $\qmin g =
g^{-1}$, we likewise get $\iptil{x}{\qmin gx} =
\iptil{gx}{x} = -\iptil{x}{gx}$. In particular
\[ \ip{\lambda_1}{\lambda_1} =
\iptil{\sigma}{(1-\qmin) \sigma} =
\iptil{\sigma}{\sigma} -
\iptil{\sigma}{\qmin\sigma} = 0-0 = 0. \]

Next, for $g=\qi$, $\qj$, or $\qk$ we have
$\ip{\lambda_1}{g\lambda_1} = \iptil{\sigma}{(1
- \qmin)g\sigma} = 2 \iptil{\sigma}{g\sigma}$;
we will compute each case separately. For
convenience, we will simplify Figure~\ref{F:
loz} to $\sigdiag{1}{\qi}{\qj}{\qmin}$, and we
will likewise represent $\qk \sigma$ by the
diagram $\sigdiag{\qk}{\qj}{\qmin\qi}{\qmin\qk}$
etc. An element $g\in G$ appears in the diagram
for a translate $g'\sigma$ of $\sigma$ if and
only if $g'\sigma$ passes through $g \vtil$.
Likewise one can reconstruct the $1$\/-cells
participating in $g'\sigma$ (with their signs).
\vspace{0.2cm}

\emph{The case $g=\qk$}\/:\ \ The four
$1$\/-cells $\altil$, $\qi \altil$, $\gatil$,
and $\qj \gatil$ in the support of $\sigma$ are
distinct from their translates by $\qk$. Hence
$\sigma$ and $\qk \sigma$ can intersect only at
points of $\Ctil$ over $v$. As is clear from the
diagrams for $\sigma$ and for $\qk\sigma$, these
points of intersection are the translates of
$\vtil$ by $\{ 1,\qi,\qj,\qmin\} \cap \{\qk,\qj,
\qmin\qi, \qmin\qk\} = \{\qj\}$. By
Figure~\ref{F: fund} the local picture at
$\qj\vtil$, when lifted to $C^{\rm univ}$ and
translated to $\hv$, is the left part of
Figure~\ref{F: cask}:
\begin{figure}[!ht]\caption{} \label{F: cask}
\setlength{\unitlength}{4mm}
\begin{picture}(31,12)(0,0)
\put(1,6){\line(1,0){10}}
\put(6,6){\vector(1,0){4}}
\put(2,6){\vector(1,0){0}}
\put(6,1){\line(0,1){10}}
\put(6,1){\vector(0,1){1}}
\put(6,9){\vector(0,1){1}}
\put(2,7){\makebox(0,0){$\hga$}}
\put(2,10.4){\makebox(0,0){$\sigma$}}
\put(5,10.5){\makebox(0,0){$\gamma^{-1} \hga$}}
\put(2,7.7){\line(0,1){2}}
\put(2.6,10.5){\line(1,0){1}}
\put(7.3,2){\makebox(0,0){$\alpha^{-1}\hal$}}
\put(10.3,2){\makebox(0,0){$\qk \sigma$}}
\put(10.3,5.5){\makebox(0,0){$\hal$}}
\put(8.7,1.9){\line(1,0){0.7}}
\put(10.3,2.8){\line(0,1){2}}
\put(14,6){\vector(1,0){4}}
\put(20,6){\line(1,0){3}}
\put(23,6){\line(1,1){2}} \put
(25,8){\line(0,1){3}} \put(25,1){\line(0,1){3}}
\put(25,4){\line(1,1){2}}
\put(27,6){\line(1,0){3}}
\put(23,6){\vector(-1,0){1.5}}
\put(25,11){\vector(0,-1){1.5}}
\put(25,1){\vector(0,1){1.5}}
\put(27,6){\vector(1,0){1.5}}
\end{picture}
\end{figure}

The right hand side of Figure~\ref{F: cask}
shows that after a homotopy $\sigma$ and $\qk
\sigma$ do not meet. Hence
$\ip{\lambda_1}{\qk\lambda_1} =
2\iptil{\sigma}{\qk \sigma} = 0$.

\emph{The case $g=\qi$}\/:\ \ Here the supports
of $\sigma$ and of $\qi \sigma \lara
\sigdiag{\qi}{\qmin}{\qk}{\qmin\qi}$ intersect
along $\qi \altil$. The left part of
Figure~\ref{F: casi} shows a neighborhood of
$\qi \altil$ lifted to $\Ctil$: the part of
$\sigma$ represented in it is $\altil +
\qi\altil - \qj\gatil$ and that of $\qi\sigma$
is $-\qi\gatil + \qi\altil + \qmin\altil$. To
obtain this picture we combine the local
pictures offered by Figure~\ref{F: fund} at both
endpoints $\qi\vtil$ and $\qmin\vtil$ of
$\qi\altil$. The right hand part of
Figure~\ref{F: casi} represents homologous
paths, and it follows that
$\ip{\lambda_1}{\qi\lambda_1} =
2\iptil{\sigma}{\qi\sigma} = -2$.

\begin{figure}[!h] \caption{} \label{F: casi}
\setlength{\unitlength}{4mm}
\begin{picture}(30,8)(0,0)
\put(1,5){\line(1,0){12}}
\put(5,5){\vector(-1,0){2}}
\put(5,5){\vector(1,0){2}}
\put(13,5){\vector(-1,0){2}}
\put(5,5){\line(0,-1){4}}
\put(5,1){\vector(0,1){2}}
\put(9,5){\line(0,-1){4}}
\put(9,5){\vector(0,-1){2}}
\put(5,1){\makebox(0,0){$\bullet$}}
\put(3,5.5){\makebox(0,0)[b]{$\qi \hga$}}
\put(7,5.5){\makebox(0,0)[b]{$\qi \hal$}}
\put(11,5.5){\makebox(0,0)[b]{$\qj \hga$}}
\put(4.5,3){\makebox(0,0)[r]{$\hal$}}
\put(9.5,3){\makebox(0,0)[l]{$\qmin \hal$}}
\put(4.5,1){\makebox(0,0)[r]{$\vtil$}}
\put(15,5){\vector(1,0){4}}
\put(21,7){\line(1,-1){4}}
\put(21,7){\vector(1,-1){3.5}}
\put(21,3){\line(1,1){4}}
\put(21,3){\vector(1,1){3.5}}
\put(21.5,3){\makebox(0,0)[l]{$\sigma$}}
\put(25.5,3){\makebox(0,0)[l]{$\qi \sigma$}}
\end{picture}
\end{figure}

\emph{The case $g=\qj$\,}\/:\ \ Here $\qj
\sigma$ corresponds to
$\sigdiag{\qj}{\qmin\qk}{\qmin}{\qmin\qj}$, so
$\sigma$ and $\qj \sigma$ intersect along $\qj
\gatil$ as in Figure~\ref{F: casj}. We handle it
exactly as we handled Figure~\ref{F: casi}, to
obtain $\ip{\lambda_1}{\qj \lambda_1} =
2\iptil{\sigma}{\qj \sigma}=-2$.
\begin{figure}[!ht]\caption{} \label{F: casj}
\setlength{\unitlength}{4mm}
\begin{picture}(30,8)(0,0)
\put(1,5){\line(1,0){12}}
\put(5,5){\vector(-1,0){2}}
\put(5,5){\vector(1,0){2}}
\put(13,5){\vector(-1,0){2}}
\put(5,5){\line(0,-1){4}}
\put(5,1){\vector(0,1){2}}
\put(9,5){\line(0,-1){4}}
\put(9,5){\vector(0,-1){2}}
\put(5,1){\makebox(0,0){$\bullet$}}
\put(3,5.5){\makebox(0,0)[b]{$\qj \hal$}}
\put(7,5.5){\makebox(0,0)[b]{$\qj \hga$}}
\put(11,5.5){\makebox(0,0)[b]{$\qi \hal$}}
\put(4.5,3){\makebox(0,0)[r]{$\hga$}}
\put(9.5,3){\makebox(0,0)[l]{$\qmin \hga$}}
\put(4.5,1){\makebox(0,0)[r]{$\vtil$}}
\put(15,5){\vector(1,0){4}}
\put(21,7){\line(1,-1){4}}
\put(25,3){\vector(-1,1){3.5}}
\put(21,3){\line(1,1){4}}
\put(25,7){\vector(-1,-1){3.5}}
\put(21.5,3){\makebox(0,0)[l]{$\sigma$}}
\put(25.5,3){\makebox(0,0)[l]{$\qj \sigma$}}
\end{picture}
\end{figure}

Now set $v_{11}=t_1 + t_2 \qi + t_3 \qj + t_4
\qk$, with $t_i\in \RR$. We get $2t_1= \Tr
v_{11}=0$, $-2t_4= \Tr v_{11} \qk= 0$, $-2t_2=
\Tr v_{11} \qi= -2$ and $-2t_3= \Tr v_{11} \qj=
-2$, so $v_{11}= \qi+\qj$.
 \vspace{0.3cm}

\emph{Computation of $v_{12}$}\/:\ \ We shall
evaluate $\iptil{\sigma}{g\hbe_{\ZG}}$ for all
$g\in G$. The loops $\sigma$ and $g\betil $ can
only intersect at points of $\Ctil$ above
$\vtil$. These intersection points are the
translates of $\vtil$ by the elements of $\{ 1,
\qi, \qj, \qmin\} \cap \{g\}$. Figure~\ref{F:
bez} shows how $\sigma$ intersects the four
translates $g \betil$, for $g\in \{1,
\qi,\qj,\qmin\}$. To verify it one translates
the basic Figure~\ref{F: fund} to the four
points above $v$ in the support of $\sigma$.
\begin{figure}[!h]\caption{} \label{F: bez}
\setlength{\unitlength}{3.5mm}
\begin{picture}(24,14)(-12,-7)
\qbezier(-10,0)(0,10)(10,0)
\qbezier(-10,0)(0,-10)(10,0)
\put(-5,3.75){\vector(2,1){0}}
\put(5,3.75){\vector(-2,1){0}}
\put(-5,-3.75){\vector(-2,1){0}}
\put(5,-3.75){\vector(2,1){0}}
\put(-2,3){\line(1,1){2}}
\put(-2,3){\vector(1,1){1}}
\put(0,5){\line(-1,1){2}}
\put(0,5){\vector(-1,1){1}}
\put(12,0){\line(-1,0){2}}
\put(12,0){\vector(-1,0){1}}
\put(10,0){\line(0,1){2}}
\put(10,0){\vector(0,1){1}}
\put(-8,0){\line(-1,0){2}}
\put(-8,0){\vector(-1,0){1}}
\put(-10,0){\line(0,-1){2}}
\put(-10,0){\vector(0,-1){1}}
\put(0,-5){\line(-1,1){2}}
\put(0,-5){\vector(-1,1){1}}
\put(2,-3){\line(-1,-1){2}}
\put(2,-3){\vector(-1,-1){1}}
\put(-5.5,3.75){\makebox(0,0)[br]{$\qi \altil$}}
\put(-1.5,6){\makebox(0,0)[tr]{$\betil$}}
\put(5.5,3.75){\makebox(0,0)[bl]{$\qj \gatil$}}
\put(-5.5,-3.75){\makebox(0,0)[tr]{$\altil$}}
\put(5.5,-3.75){\makebox(0,0)[tl]{$\gatil$}}
\put(0,-4.5){\makebox(0,0)[b]{$\betil$}}
\put(-8.5,-0.1){\makebox(0,0)[tl]{$\qi \betil$}}
\put(10.5,0.5){\makebox(0,0)[bl]{$\qj \betil$}}
\put(0,-5){\makebox(0,0){$\bullet$}}
\put(0,-5.5){\makebox(0,0)[t]{$\vtil$}}
\end{picture}
\end{figure}

This implies that
$\iptil{\sigma}{\hbe_{\ZG}}=0$,
$\iptil{\sigma}{\qi\hbe_{\ZG}}=1$,
$\iptil{\sigma}{\qj\hbe_{\ZG}}=0$, and
$\iptil{\sigma}{\qmin\hbe_{\ZG}}=1$. The other
intersections $\iptil{\sigma}{g\hbe_{\ZG}}$ are
trivially $0$. Hence for $g=1$, $\qi$, $\qj$,
and $\qk$ we have $\ip{\lambda_1}{g\lambda_2} =
\ip{\sigma}{g\beta} - \ip{\sigma}{\qmin g \beta}
= -1$, $1$, $0$, and $0$ respectively. Writing
$v_{12}=t_1+t_2\qi +t_3 \qj+t_4 \qk$ we see as
before that $v_{12}=(-1-i)/2$. Hence
$v_{21}=(1-i)/2$. This concludes the proof of
Theorem~\ref{PEL}.
\end{proof}

\begin{remark}
\label{standard} Set $\lambda'_1 =
-(1+\qi)\lambda_1 + (\qi+\qk)\lambda_2$ and
$\lambda'_2 = \lambda_2$. Then $\lambda'_1$ and
$\lambda'_2$ constitute a $B$-basis for
$H_1(P,\QQ)$. The pairing is given in this basis
as in \eqref{E: pair} but with the matrix
$[v'_{ij}]= \matr{0}{1}{-1}{0}$.

We omit the routine verification.
\end{remark}

The ideal $\bP = (1+\qi)\bM$ is two-sided, of
index $4$ in $\bM$, and contained in $\bM'$ with
index $2$. Let $\ip{\;}{\;}_{\bM'}: \bM' \times
\bM' \ra \ZZ$ map $(m_1,m_2)$ to $\frac{1}{2}
\Tr_{B/\QQ}( \ov{m_1} m_2)$. We will prove the
following:
\begin{theorem}
\label{PEL1} For an elliptic curve $E/\CC$ let
$A_E$ be the $4$\/-dimensional abelian variety
$\bM' \otimes E$ with the polarization
\[ \ip{\;}{\;} = \ip{\;}{\;}_{\bM'} \otimes
\ip{\;}{\;}_E ,\] where $\ip{\;}{\;}_E$ is the
standard polarization on $E$. Then
\begin{enumerate}
\item The polarization $\ip{\;}{\;}$ is
principal. It is preserved by $\bM'$ up to
similitudes. As a polarized variety $A_E$ is
isomorphic to $E^4$  with the product
polarization.
\item Let $\alpha: E\ra E'$ be an isogeny of
degree $2$ of elliptic curves, so that
$H_1(E',\ZZ)$ contains $H_1(E,\ZZ)$ with index
$2$. In the isogeny class of $A_E$, let
$B_\alpha$ be the polarized abelian variety
characterized by
\begin{equation}
\label{E:Bal}
 H_1(B_\alpha,\ZZ) = 2\bM \otimes H_1(E',\ZZ) +
 \bP \otimes H_1(E,\ZZ) \subset H_1(A,\QQ).
\end{equation}
Then the induced pairing on $H_1(B_\alpha,\ZZ)$
makes $B_\alpha$ into a principally polarized
abelian variety with an $\bM$\/-action.
\item Each Prym variety $\Prym(\Ctil/C_\pm)$
in case $\cM(2,0,0,0)$ is isomorphic to some
$B_\alpha$ as a polarized abelian variety
together with the $\bM$\/-action.
\end{enumerate}
\end{theorem}
\begin{proof}
The first part is clear when using the standard
basis $1$, $\qi$, $\qj$ and $\qk$ of $\bM'$,
since $\ip{\;}{\;}_{\bM'}$ is then the standard
Euclidean pairing (and $\bM'$ acts via
similitudes). For the second part, let $e$, $f$
be a symplectic basis for $H_1(E,\ZZ)$ so that
$\frac{1}{2}e$ and $f$ are a basis for
$H_1(E',\ZZ)$. Then
\[
 H_1(B_\alpha,\ZZ) = \bM\otimes e \oplus
\bP\otimes f \subset H_1(A,\QQ)\otimes B.
\]
Observe that $\Tr_{B/\QQ}$ is even on $\bP$. It
follows that $\ip{\;}{\;}$ is integral on
$H_1(B_\alpha,\ZZ)$, hence it defines a
polarization on $B_\alpha$. Since the lattices
$H_1(B_\alpha,\ZZ)$ and $H_1(A_E,\ZZ)$ have the
same volume, it follows that this polarization
is principal as asserted. The $\bM$\/-action
preserves $H_1(B_\alpha,\ZZ)$ since $\bP$ is an
$\bM$\/-ideal, and the second part follows.

To prove the third part, we use
Remark~\ref{standard} to write $H_1(P,\QQ) =
B\otimes L$, where $L = \ZZ\lambda'_1\oplus
\ZZ\lambda'_2$. By a standard exceptional
isomorphism of Lie groups, the group $J$ of
$B_\RR$\/-linear similitudes of $H_1(A_E,\RR)$
is then $\GL(L\otimes\RR)\times B^\ast/\sim$,
where we identify $(t,1)\sim(1,t)$ for any
scalar $t$ (see e. g. \cite[Chap.
IX.4.B.xi]{Helg}). Let $h_P:\CC^\ast \ra
J\subset \GSp(H_1(P,\RR),\ip{\;}{\;})$ be the
Hodge type of $P$, in the sense of \cite[Section
4]{del}. Since $h_P(\sqrt{-1}\,)$ is a Cartan
involution of $J$, the image of $h_P$ must
centralize the compact factor subgroup of $J$
consisting of the norm $1$ elements in $B^\ast$.
Therefore $h_P$ factorizes through a Hodge type
$h_0 = h_P:\CC^\ast \ra \GL(L\otimes\RR)$. The
lattice $L\subset L\otimes \RR$ then determines
an elliptic curve, characterized by the
properties that its Hodge type is $h_0$, and
that $H_1(E,\ZZ) = L$. The over-lattice $L' =
\ZZ\lambda_1\oplus \ZZ\lambda_2$  defines
similarly an elliptic curve $E'$ with an isogeny
$\alpha:E'\ra E$ of degree $2$. From the
definitions, the resulting equality
$H_1(\Prym,\QQ) = B\otimes L$ is compatible with
the $B$\/-action, with the Hodge structure, and
with the ($B$\/-hermitian) polarization. Now
define $B_\alpha$ by~\eqref{E:Bal} above. Since
$\lambda'_1$ and $\lambda'_2$ are a symplectic
basis for $L$, it follows that we have an
equality of lattices in this space
\[
H_1(\Prym,\ZZ) = \bM\lambda_1 \oplus
\bM\lambda_2 = \bP\lambda'_1 \oplus
\lambda'_2\bM = H_1(B_\alpha,\ZZ),
\]
so that $\Prym \simeq B_\alpha$ as asserted.
\end{proof}
\begin{corollary}
\label{PEL2} Let $Y_0(2)$ be the modular curve
parameterizing elliptic curves with an isogeny
$\alpha:E\ra E'$ of degree $2$, and let $w_2$ be
the modular (Atkin-Lehner) involution of
$Y_0(2)$, sending $\alpha$ to its dual isogeny.
Then the quotient curve $Y_0(2)/w_2$ is
isomorphic to the Shimura curve $\Shim$
parameterizing the PEL data of
Theorem~\ref{PEL1}.2 (see~\cite{del, shi}) via
the assignment $\phi: \alpha\mapsto B_\alpha$.
\end{corollary}
\begin{proof}
By Theorem~\ref{PEL1}(2), $\phi:Y_0(2) \ra
\Shim$ is a morphism. Moreover by
Theorem~\ref{PEL1}(3) $\phi$ is surjective.
Since we work over $\CC$ we know that
analytically $Y_0(2) = \Gamma_0(2)\backslash
\cH$, where $\Gamma_0(2)$ is the subgroup of
$\SL(2,\ZZ)$ consisting of the matrices
$\matr{a}{b}{c}{d}\in\SL(2,\ZZ)$ for which
$2|c$. Since $\phi$ is modular it induces an
isomorphism $\Gamma\backslash\cH\sra{\sim}\Shim$
for some congruence subgroup $\Gamma$ of
$\PSL(2,\RR)$ containing $\Gamma_0(2)$. Let $e$,
$f$ be the symplectic basis for $L$ as in the
proof of Theorem~\ref{PEL1}(2). In terms of this
basis Let $W_2: B\otimes L\ra B\otimes L$ be the
involution $R_{(1+\qi)/2} \otimes
\matr{0}{-1}{2}{0}$, where $R_{(1+\qi)/2}$ acts
as right multiplication by $(1+\qi)/2$. Then
$W_2$ is clearly (left) $B$\/-linear. Moreover,
it preserves $H_1(B_\alpha,\ZZ)$ and the
polarization, and its $L$ component generates
the normalizer $N_0(2)$ of $\Gamma_0(2)\subset
\SL(L)$. To see this, we compute
\[W_2(\bM e \oplus \bP f) = \bP(1+\qi)/2 e
\oplus 2\bM(1+\qi) f = \bM e \oplus \bP f; \quad
\text{likewise,}
\]
\[\begin{array}{rl}
\ip{W_2(m_1e+m_2f)}{W_2(n_1e+n_2f)} &=
\ip{m_2\frac{1+\qi}{2}e + m_1(1+\qi)f}
{-n_2\frac{1+\qi}{2}e + n_1(1+\qi)f} \\
&=\Tr_{B/\QQ}
(-\frac{1-\qi}{2}\ov{m_2}n_1(1+\qi) +
(1-\qi)\ov{m_1}n_2\frac{1+\qi}{2})\\
&=\Tr_{B/\QQ}(\ov{m_1}n_2-\ov{m_2}n_1)\\
&=\ip{m_1e+m_2f}{n_1e+n_2f}.
\end{array}\]
The last part is well-known. Hence our $\Gamma$
contains $N_0(2)$. Since $N_0(2)/\pm1$ is known
to be maximal as a fuchsian group, we get
$\Gamma = N_0(2)$, proving our assertion.
\end{proof}
\begin{remark}
\label{GV2} Parts 1 and 2 of Theorem~\ref{PEL1}
are stated in \cite[Proposition 2.6]{GV}, and
their approach is a geometric version of our
explicit argument for Theorem~\ref{PEL1}.3.
However, the lengthy analysis which we needed to
determine the pairing and deduce its properties,
which in their terminology would have amounted
to the analysis of the contraction map, is not
done in their paper.
\end{remark}

\section{Cubics with nine nodes} \label{cubics}
In this section we will study cubic threefolds
$X$ with nine nodes (i. e. ordinary double
points). In the next section these will be
related to the quaternionic abelian varieties
through some Prym-theoretic constructions.

The maximal number of nodes that a cubic
threefold $X$ can have is $10$, and this happens
if and only if $X$ is (projectively) the Segre
cubic (see \cite{seg, var, don} and
Lemma~\ref{equ} below). We thank Igor Dolgachev
for telling us about the beautiful work
\cite{seg}. In it, C. Segre studies cubics with
$n$ nodes, $6 \leq n \leq 10$. He starts with
the subvarieties $S \subset D \subset \PP^8$,
where $\PP^8$ is the projectivization of the
vector space of 3x3 matrices, $D$ is the locus
$\mbox{det}=0$ of singular matrices, and $S$ is
the locus of rank-1 matrices which we would
nowadays call the Segre embedding $\PP^2 \times
\PP^2 \hookrightarrow \PP^8$. The intersection
of $D$ with a generic subspace $\PP^4 \subset
\PP^8$ is a cubic threefold with 6 nodes at the
points of $\PP^4 \cap S$. By moving the $\PP^4$
subspace into special position, the cubic
threefold $\PP^4 \cap D$ can be made to have $n$
nodes, $6 \leq n \leq 10$. Segre states that
conversely, any cubic with $n \geq 6$ isolated
nodes can be obtained this way, and proceeds to
a detailed case by case analysis. In the case of
interest to us, he points out that 9-nodal
cubics appear in the pencil generated by two
completely reducible cubics, i.e. $x_1 x_2 x_3 +
\alpha x_4 x_5 x_6=0$, where the $x_i$ are six
general linear coordinates on $\PP^4$ satisfying
a single linear relation which we can write as
$\sum_i x_i =0$. He states that all 9-nodal
cubics arise this way. The properties of such
cubics are then straightforward to determine.

In this section we give a modern treatment of
these results and explain their modular
interpretation. The cubics with $9$ nodes turn
out to form an irreducible family with many nice
properties. In fact we have the following:
\begin{theorem}
\label{nine} {\bf (1)\ } Let $X$ be a cubic
threefold with at least $9$ {\em isolated}
singularities over an algebraically closed
field. Then the singular locus $X_\sing$ of $X$
consists of \ $9$ or $10$ nodes.

\noindent {\bf (2)\ } Through every node of an
$X$ as in (1)\ pass $4$ planes contained in $X$.
Each plane $P'$ contained in a cubic threefold
$X'$ with at most isolated singularities along
$P'$ contains at most $4$ singularities of $X$.
It contains exactly $4$ singularities if and
only if they are all nodes for $X$.

\end{theorem}
\begin{proof}
Fix an isolated singularity $O$ of $X$. The
lines through $O$ contained in $X$ form an
algebraic set $C_{O,X}$ in the projectivized
tangent space $\PP = \PP(T_O(\PP^4)) \simeq
\PP^3$. Choose homogeneous coordinates
$[x;y;z;w;u]$ of $\PP^4$, where $O =
[0;0;0;0;1]$. Then $X$ is given by an equation
$f = uq + c = 0$, and $C_{O,X}$ by $q = c =0$,
where $q$ and $c$ are a quadric and a cubic in
$(x,y,z,w)$ respectively. To relate the
singularities of $X$ and of $C_{O,X}$, let
$\nabla$ denote the gradient in the $x,y,z,w$
variables. We will need the following facts:
\begin{lemma}
\label{circast}
 Assume that $X$ is given by $f=uq+c$
as above, with $O = [0;0;0;0;1]$ an isolated
singularity as before. Then

\noindent{\bf (1)\ } The singularities of $X$
are given by $q=c=0$ and $u\nabla q+\nabla c =
0$. The singularities of $C_{O,X}$ are the
points of $\PP$ where $\nabla q$ and $\nabla c$
are linearly dependent.

\noindent{\bf (2)\ }  Let $P\neq O$ be any
singularity of $X$, and let $p\in\PP$ denote the
point corresponding to the line $\ov{OP}$. Then
$q$ is nonsingular at $p$, i. e. $\nabla q(p)
\neq 0$.
\end{lemma}
\begin{proof}
The first part is immediate. For the second
part, put $p$ at $[0,0,0,1,0]$. If $q$ were
singular at $p$ then $q=q(x,y,z)$ would not
depend on $w$. Moreover we would have
$c=c'(x,y,z)+wq'(x,y,z)$. But then $X$ would be
singular along the entire line $\ov{OP}$,
contradicting our assumption.
\end{proof}
Note that $C_{O,X}$ is a curve: otherwise $q$
and $c$ would have a component in common,
yielding $f=L q'$ for a linear form $L$ and a
quadric $q'$. But then $X$ would have at most
one isolated singularity (the vertex of the cone
$q'=0$\/). Observe also that a line $\ell$
through two singularities of any cubic threefold
$Y$ given by $\{g=0\}$ is contained in $Y$, and
that a line through three singularities of $Y$
is contained in $Y_\sing$: indeed, a cubic
polynomial $g$ defining $Y$ vanishes in the
first case to order $2$ at the $2$ singularities
so it vanishes on $\ell$, while in the second
case each partial derivative of $g$ vanishes at
three points of $\ell$ and has degree $2$, so is
$\equiv 0$.

Returning to our $X$ we see from this and from
Lemma~\ref{circast} that a line on $X$ through
$O$ corresponds to an (isolated) singularity of
$C_{O,X}$ if and only if it contains an
(isolated) singularity
--- necessarily unique --- of $X$, which is
different from $O$. In particular $X$ has $m+1$
isolated singularities if $C_{O,X}$ has $m>0$
isolated singularities.

To prove part 1 it now suffices to prove that
$q$ has maximal rank $r=4$, so that $O$ is a
node,  and that $C_{O,X}$ has at most one
additional singularity to the $8$ we know.

Notice first that $r\geq 2$ by
Lemma~\ref{circast}(2). Suppose next that $r =
2$. Then we may write $q = wz$. The
singularities of $X$ are on the union of the two
hyperplanes $z =0$ and $w = 0$ but not on their
intersection (by Lemma~\ref{circast}(2)). Thus
$C_{O,X}$ is contained in the union of the
corresponding planes in $\PP$. On each plane
$C_{O,X}$ is defined by the cubic equation $c
=0$, so it can have at most $3$ isolated
singularities. Since $C_{O,X}$ cannot be
singular along the intersection $z=w=0$, we see
that altogether $C_{O,X}$ has at most $6$
isolated singularities, contradicting our
assumption.

Suppose now that $r=3$. Then $q=0$ defines a
quadric cone $S_0$ in $\PP$, whose vertex $e$
cannot be in $(C_{O,X})_\sing$ by
Lemma~\ref{circast}(2). We shall obtain a
contradiction by showing that the intersection
$C'$ of a cubic surface in $\PP$ with $S_0$ can
have at most $6$ (isolated) singularities away
from $e$. For this let $S$ be the blowup of
$S_0$ at $e$. It is well-known that $\Pic{S}$ is
freely generated by the exceptional divisor $E$
and the proper transform $F$ of a line on $S$
through $e$. Let $H$ be the pullback to $S$ of
the hyperplane class $\cO_\PP(1)_{|S}$. Then $2F
+E \equiv H$ in $\Pic{S}$. A divisor class $aE +
bF$ contains a reduced and irreducible curve $B$
if and only if either $(a,b)=(0,1)$ (and then
$B$ is a fiber $F$\/), or $(a,b) = (1,0)$ (and
then $B$ is $E$\/), or if $b\geq 2a>0$ (and then
$B$ is in $|aH+cF|$ with $c\geq 0$\/). Indeed
intersecting $B$ with $E$ and with $F$ shows
these conditions are necessary, and their
sufficiency follows from Bertini's theorem,
since the general member in $aH+cF$ is smooth,
hence irreducible. The canonical class is $K_S =
-2E-4F = -2H$ and hence $K_S^2 = 8$. By
adjunction, the arithmetic genus of an
irreducible curve in $|aE+bF|$ as above is $0$
in the first two cases $(a,b)=(1,0),(0,1)$, is
$0$ if $a=1$ and $c\geq 0$, is $1+c$ if $a=2$
and $c\geq 0$, and is $4$ if $(a,b)=(3,6)$.
Since an irreducible curve of arithmetic genus
$g$ has at most $g$ singularities, a member of
$H+cF$ is smooth.

We will now show that the proper transform $C''$
of $C'$ has at most $6$ nodes not on $E$ by
examining the types of irreducible components
that $C''\in|3H|$ can have. Each component is
$E$, or some $F$, or of type $aH+cF$ with $3\geq
a>0$. A component with $a=3$  can have at most
$4$ singularities, a component with $a=2$ has at
most $1$ singularity, and a component with $a=0$
is nonsingular. In particular, since a component
with $a=3$ must be all $C''$ and $4<7$, this
case cannot occur. Similarly, a component $D$
with $a=2$ cannot occur, since the other
components are either $H$
--- then $D\in|2H|$, and there are at most one
singularity on $D$, none on $H$ and four points
of intersection, making a total of at most $5<7$
nodes. Otherwise, there are $2-c$ components of
type $F$ and one component equal to $E$, and
there are at most $1+2(2-c)<7$ nodes. In
conclusion, only $a=1$ occurs. Next, the number
$k$ of components of type $E$ clearly cannot be
more than $3$. It cannot be $3$ since the other
components will be only $F$\/'s and $E$, without
any nodes not on $E$; $k$ cannot be $2$, since
then we will have one component of type $H+cF$
and $4-k$ fibers, giving at most $4-k<7$ nodes;
if $k=1$ we have components $D_1$, $D_2$ of
types $aH+c_iF$, $i=1$, $2$, and $2-c_1-c_2$
fibers: this gives at most $5<7$ nodes. Finally,
when $k=0$ we get $3$ components of type $H$ and
there are at most $6$ nodes.

We now know that $r=4$, so the locus of $q=0$ is
a nonsingular quadric $S$. Let $F$, $F'$ be the
two standard rulings of $S$ by lines. As before
we want to find the maximal number of isolated
singularities $p_i$ that a member $C'$ of
$3(F+F')$ can have. Applying the same type of
analysis as before we find that when this number
is eight or more, $C'$ breaks into two fibers of
$F$, two fibers of $F'$, and a member of the
hyperplane class $H \simeq F + F'$, all
intersecting transversely. Moreover $H$ is
reducible if and only if $C'$ has $9$ nodes.
Part 1.\ of the Theorem follows.

The explicit description of such a curve $C'$ in
our case $C'= C_{O,X}$ shows that on each of the
$4$ line components $l$ of $C_{O,X}$ there are
$3$ singularities. Hence each of the planes
$\ov{Ol}$ contains the three corresponding nodes
of $X$ in addition to the node $O$. This plane
intersects $X$ in at least the $6$ lines joining
any two of these $4$ nodes of $X$, hence is
contained in $X$. This gives $4$ planes through
$O$ contained in $X$.

Finally, let $\Pi'$ be a plane contained in a
cubic threefold $Y:\{g=0\}$. If the plane is
given by $w=z=0$, then the singularities of $Y$
along $\Pi'$ are the intersection of the two
conics $\partial g/\partial z = \partial
g/\partial w = 0$. A point in the intersection
is a node for $Y$ if and only if the
intersection is transverse there,  and all the
intersections are transverse if and only if
there are precisely $4$ of them, proving part 2
of the Theorem.
\end{proof}

Using the Theorem we can describe the cubic
threefolds having at least $9$ nodes:
\begin{lemma}
\label{equ} A cubic threefold over a scheme $S$
containing nine given nodes is projectively
equivalent to one given by
\begin{equation}
\label{E:nine} X(\alpha):\qquad x_1 x_2 x_3 +
\alpha x_4 x_5 x_6 = x_1+ \dots +x_6 = 0
\end{equation}
in $\PP^5/S$, where $\alpha$ is in $\GG_m(S)$.
Under the evident $S_3\times S_3$ symmetry, the
nine nodes are the orbit of $O_{3,6} =
(0,0,1,0,0,-1)$. Over an algebraically closed
field there are $10$ nodes precisely in the
Segre case $a=1$, and then the $10$\/th node is
$(1,1,1,-1,-1,-1)$.
\end{lemma}
\begin{proof}
Take affine coordinates $x'_1,x'_2,x'_4,x'_5$
for $\Aa^4$ so that the origin $O$ is one of the
given nodes of $X$. Let $T/S$ be the locus in
the grassmanian $\Gr(2,T_O(\Aa^4))/S$ of planes
through $O$ which are contained in $X$ and whose
intersection with the singular locus $X_{{\rm
sing}}$ of $X$ is supported on the given nodes.
The explicit description of $C_{O,X}$ obtained
in the proof of Theorem~\ref{nine} gives that
$T$ is \'etale of degree $4$ over $S$. We also
know that $T$ corresponds to $4$ lines on
$C_{O,X}$, which intersect mutually according to
the graph of the sides of a square. These
intersections represent $4$ of the given nodes
of $X$, so that monodromy acts trivially on the
square, and it follows that $T$ is a trivial
(product) covering of $S$. The projectivized
tangent cone to $X$ at $O$, which is a
nonsingular quadric $Q$, contains these two
lines in each of its rulings which are marked,
i.e.. the \'etale cover of $S$ which these lines
define is a trivial (product) cover. Hence we
may choose the coordinates so that these lines
are $x'_i = x'_j=0$ for $1\leq i,j \leq 2$. The
tangent cone is then $x'_1x'_2 + tx'_4x'_5 = 0$
for some $t\in \cO_S^\times$, and replacing
$x'_4$ by $tx'_4$ we may assume $t=1$. Then
$C_{O,X}$ is the intersection of the tangent
cone above with $x'_1x'_2\ell_1=0$, where
$\ell_i$ denotes an $\cO_S$\/-linear function of
$x'_1$, $x'_2$, $x'_4$, and $x'_5$ for any $i$.
Homogenizing, we see that an equation of $X$ in
$\PP^4$ is given by
\[(x'_3+\ell_2)(x'_1x'_2+x'_4x'_5) +
x'_1x'_2\ell_1 = 0.
\]
Taking $y_3 = x'_3 + \ell_2 + \ell_1$, $y_i =
x'_i$ for $i=1$, $2$, $4$, and $5$, and $m =
x'_3 + \ell_2$ gives the equation $y_1 y_2 y_3 +
y_4 y_5 \sum_{i=1}^5 \alpha_i y_i=0$. The
coefficients $\alpha_i\in\cO_S$ are invertible
on $S$: if $\alpha_1(s) =0$ for a geometric
point $s$ of $S$, then the tangent cone to the
singularity $O_3=[0;0;1;0;0]$ is reducible, and
similarly for $i=2$ or $3$; if $\alpha_4=0$
(respectively $\alpha_5=0$ then $O_4 =
[0;0;0;1;0]$ (respectively $O_5=[0;0;0;0;1]$) is
a singularity whose tangent cone is reducible.
In each case we get a non-nodal singularity on
$X$, contradicting our assumption. Thus
$\alpha_i$ is invertible, so we may replace each
$y_i$ by $x_i = \alpha_i y_i$. Setting $x_6 =
=-\sum_{i=1}^5 x_i$ we get the desired form.

To determine the singularities we must find the
points when the gradients of the two equations
in~\eqref{E:nine} defining $X(\alpha)$ are
dependent. If any $x_i$ is $0$ we get that two
of $x_1$, $x_2$, $x_3$, and two of $x_4$, $x_5$,
and $x_6$ are $0$. This leads to the nine nodes
of type $O_{3,6}$. Else we find the 10th node as
indicated with $\alpha=1$. We omit the details.
\end{proof}

Note that our formulas are characteristic free.
In characteristic $>3$ the Segre cubic threefold
is usually given by the $S_6$\/-symmetric
equations
\[\sum_{i=1}^6 y_i = \sum_{i=1}^6 y_i^3 = 0.\]
The coordinate change $y_i=x_j+x_k-x_i$, for
$\{i,j,k\} = \{1,2,3\}$ or $\{4,5,6\}$,
transforms our form into the other in an
$(S_3\times S_3)\rtimes S_2$\/-equivariant way.

We will show that the function $\alpha$ in
Lemma~\ref{equ} is unique and that the
family~\eqref{nine} is universal; this will
require (a part of) the following
\begin{proposition}
\label{PS} Let $X$ be a cubic threefold with
nine given nodes (over any base, with some
fibers possibly having a $10$\/th node). Then we
have the following:
\begin{enumerate}
\item $X$ contains $9$ planes in bijection with
the nodes, with a node $p$ corresponding to a
plane $\Pi$ if for every other plane $\Pi'$ we
have $p\in\Pi'\Leftrightarrow \Pi'$ is
transversal to $\Pi$. In particular, if $X$ has
exactly $9$ nodes, then it contains exactly
(these) $9$ planes.
\item There are exactly six {\em plane systems},
namely sets of three pairwise transverse planes
of these nine on $X$.
\item Two plane different systems are disjoint
or have one plane in common.
\item If we define two plane systems to be
equivalent whenever they are equal or disjoint,
then this is indeed an equivalence relation, and
there are two equivalence classes $A,B$
consisting of three plane systems each. In this
way the planes in $X$ are put in bijection with
$A\times B$: each plane is in a unique system of
type $A$, and in a unique system of type $B$.
\end{enumerate}
\end{proposition}
\begin{proof}
In the coordinates of Lemma~\ref{equ}, the node
$O_{i,j}$ having $1$ at the $i$'th coordinate,
$-1$ at the $j$'th coordinate, and $0$ elsewhere
corresponds to the plane $\Pi_{i,j} = \{x_i =
x_j = 0\}$ for any $1\leq i\leq 3$ and $4\leq
j\leq 6$. The plane systems of class $A$ consist
of the planes $x_i = x_j = 0$ with $1\leq i\leq
3$ fixed and each of $4\leq j\leq 6$; those of
class $B$ consist of the planes $x_i = x_j = 0$
with each of $1\leq i\leq 3$ and $4\leq j\leq 6$
fixed. All the assertions of our Proposition are
now straightforward.
\end{proof}

\begin{definition}
Set $A=\{1,2,3\}$ and $B=\{4,5,6\}$. An {\em
allowable marking} of a cubic threefold with
nine given nodes is a bijection of the given
nodes with $A\times B$ for which there exist the
nodes-planes configuration indexed as in
Proposition~\ref{PS}. For $a\in A$ and $b\in B$
we will mark the corresponding node by $O_{ab}$.
\end{definition}
We can now strengthen Lemma~\ref{equ}:
\begin{theorem}
\label{moduli} \ \noindent {\bf (1)\ } The
moduli problem of classifying nine-nodal cubic
threefolds with an allowable marking is
represented by $\GG_m$. Let $\alpha$ be the
usual coordinate of $\GG_m$. Then a universal
family is given by \eqref{E:nine}. The universal
family is (allowably) marked by letting each
marked node $O_{ab}$ of $X$ be the point of
$\PP^4$ having $1$ at the $i$\/th coordinate,
$-1$ at the $j$\/th coordinate, and $0$
elsewhere.

\noindent {\bf (2)\ } The moduli problem of
classifying cubic threefolds with nine unmarked
nodes but with marked plane systems $A$, $B$ is
coarsely represented by the same (under the
forgetting functor) $\GG_m$.

{\bf (3)\ } Over an algebraically closed field,
$X(\alpha)$ and $X(\beta)$ of \eqref{E:nine} are
isomorphic if and only if $\beta = \alpha^{\pm
1}$. The involution interchanging $A$ and $B$ on
$\GG_m$ above is given by $\alpha \mapsto
\alpha^{-1}$.

\noindent {\bf (4)\ } The moduli problem of
classifying cubic threefolds with unmarked nine
nodes is coarsely represented by
$\GG_m/(\alpha\sim \alpha^{-1}) = \Aa^1$, with
coordinate $b =\alpha + \alpha^{-1}$.
\end{theorem}
\begin{proof}
{\bf (1)\ } Let $X/S$ be a family of cubic
threefolds with $9$ nodes, marked $\{O_{a,b}\}$
as above, over a base scheme $S$. We view the
ambient $\PP^4/S$ as the hyperplane $x_1+\dots
+x_6=0$ in $\PP^5/S$. We will show that there
are unique coordinates on $\PP^4$ so that each
node $O_{ab}$ and each plane $\Pi_{ab}$ on $X$
go to its namesake in $\PP^4$: indeed, the proof
of Lemma~\ref{equ} started by doing this for
$O_{3,6}$. Then, perhaps after permuting $x_1$
with $x_2$ and/or $x_4$ with $x_5$, we got it
also for the $\Pi_{ab}$\/'s with $a=1,2$ and
$b=4,5$. As there is no pair of permutations of
$A$ and $B$ fixing these, the rigidity of the
configuration of nodes and planes on $X$ of
Proposition~\ref{PS} now forces each node and
plane of $X$ to go to its namesake in $\PP^4$ as
asserted. The $9$ nodes are in general position
in $\PP^4$, in the sense that a linear
automorphism of $\PP^4$ fixing them (pointwise)
is the identity. Hence the coordinates are
indeed unique. In other words, $X/S$ is the
pull-back of the family $X(\alpha)/\GG_m$ via a
unique morphism $S\ra\GG_m$ compatible with the
markings. This is what we had to show.

\noindent {\bf (2)\ } The $S_3\times S_3$ action
on $A\times B$ acts trivially on
$\alpha\in\GG_m$ (from part (1)), and dividing
this $\GG_m$ by the trivial action(!) gives the
claim.

\noindent {\bf (3)\ } Notice that the
automorphisms of $A\times B$ of the form
$\sigma_A\times \sigma_B$ are realized by linear
automorphisms of $\PP^4$ preserving each
$X(\alpha)$. In addition, exchanging $x_i$ with
$x_{3+i}$, for $1\leq i\leq 3$ gives an
isomorphism $\theta$ of $X(\alpha)$ with
$X(\alpha^{-1})$. This isomorphism interchanges
the classes $A$ and $B$. Now let $\phi:X(\alpha)
\ra X(\beta)$ be an isomorphism. Since the
singularities of these threefolds are in
codimension $3$, the weak Lefschetz theorem
tells us that their Picard groups are those of
the ambient projective space, namely $\ZZ$, with
the hyperplane class as canonical generator.
Hence $\phi$ must preserve it, and so is induced
by a linear automorphism of $\PP^4$. The
rigidity of the plane systems of
Proposition~\ref{PS} shows that after
composition with $\sigma_1\times\sigma_2$ and
possibly with $\theta$, our $\phi$ must map each
node $O_{ab}$ of $X(\alpha)$ to its namesake in
$X(\beta)$. As was already remarked, this forces
$\phi$ to be the identity, and in particular
$\beta = \alpha^{\pm 1}$ as asserted.

\noindent {\bf (4)\ } This follows again by
dividing $\GG_m$ by $\alpha \sim \alpha^{-1}$.
\end{proof}

Recall that the lines on a cubic threefold $X$
having at most isolated singularities form a
surface $F(X)$, called the Fano surface of $X$.
Let $O$ be a node of $X$ and as before, let
$C_{O,X}$ be the curve of lines in $X$ through
$O$ as before. If $X$ is generic (among cubics
with $O$ as a node) then $F(X)$ is identified
with $\Sym^2C_{O,X}$, where a line $\ell$ on $X$
not passing through $O$ is mapped to the two
lines through $O$ which the plane $\ov{O\ell}$
cuts on $X$. For cubic threefolds with $9$ or
$10$ nodes $F(X)$ is reducible and is described
as follows:

\begin{proposition}\label{Fanos}
Let $X$ be any cubic threefold with nine nodes
and let $X'$ be the Segre cubic threefold. Then
we have the following
\begin{enumerate}
\item $F(X')$ consists of fifteen dual planes
$\Pisp_{ij}$, for $1\leq i<j\leq 6$, and of six
{\em rulings} $R'_i$, $1\leq i\leq 6$, namely
the set of lines on $X$ meeting each plane in
the $i$\/th plane system.
\item For any plane system on $X$ let $R_i$, $1\leq
i \leq 6$ be the corresponding ruling. Then
$R_i$ is a surface.
\item $F(X)$ consists of the nine dual planes
$\Pis_{ab}$, for $(a,b)\in A\times B$, and the
six rulings $R_i$, for $i\in A\cup
B=\{1,\dots,6\}.$ If we view $F(X)$ as a cycle
on the grassmanian $G(2,\PP^4)$, then each
component counts with multiplicity $1$.
\item Under the degeneration $a\ra 1$ of $X$
to $X'$ given in Lemma~\ref{equ}, the plane
$\Pis_{ij}$ goes to the plane $\Pisp_{ij}$ for
$1\leq i \leq 3$ and $4\leq j \leq 6$, and the
ruling $R_i$ degenerates to $R'_i + \Pisp_{jk}$
whenever $\{i,j,k\} = \{1,2,3\}$ or $\{4,5,6\}$.
\end{enumerate}
\end{proposition}
\begin{proof}
The first part is well-known (see e.g.
\cite{don}). It is convenient to fix the node
$O=O_{3,6}$ and to let $C_{O,X'}$ be the curve
of lines on $X'$ through $O$ as before. Write
$C_{O,X'}$ as a union of lines $C_{O,X'}=
\cup_{i=1}^6 L_i$, where the first three are of
type $(1,0)$ and the last three of type $(0,1)$
on the quadric. Then the components of
$\Sym^2C_{O,X'}$ correspond to those of $F(X')$
by $\{L_i,L_j\} \lara \Pis_{ij}$ and $\Sym^2
L_i\lara R'_i$.

The second part follows easily from
equation~\eqref{E:nine}. For the remaining
parts, consider a degeneration of $X$ to $X'$ in
which the node $O=O_{36}$ is fixed, and where
the curve $C_{O,X}$ acquires in the limit
$C_{O,X'}$ another node. We can write $C_{O,X}=
D\cup L_{16}\cup L_{26}\cup L_{35}\cup L_{3,4}$,
where the $L$\/'s are fixed lines of types
$(1,0),(1,0),(0,1),(0,1)$ on the quadric and $D$
is the conic (of type $(1,1)$\/) degenerating to
two lines. Then the components of
$\Sym^2C_{O,X}$ correspond to those of $F(X)$ as
follows.

\begin{equation}
\label{E:sym2}
\begin{array}{l}
\mbox{ \noindent $\Sym^2 L_{ij} \lara \Pis_{ij}$
for $\{i,j\}=\{1,6\},\{2,6\},\{34\},\{3,5\}$;}
\\ \mbox{
\noindent $\{L_{i6},L_{3j}\} \lara \Pis_{k\ell}$
for $\{i,k\}=\{1,2\}$ and $\{j,\ell\}=\{4,5\}$;}
\\ \mbox{
\noindent $\{L_{34},L_{35}\} \lara R_3$ and
$\{L_{16},L_{26}\} \lara R_6$;}
\\ \mbox{
\noindent $\{L_{i6},D\} \lara R_i$ and
$\{L_{3j},D\} \lara R_j$ for $i\in \{1,2\}$ and
$j\in \{4,5\}$;}
\\ \mbox{
\noindent $\Sym^2 D \lara \Pis_{36}$. }
\end{array}
\end{equation}

Since we know how $C_{O,X}$ degenerates to
$C_{O,X'}$, we know how $F(X) = \Sym^2C_{O,X}$
degenerates to $F(X')=\Sym^2C_{O,X'}$:
$\Pis_{ij}$ is constant for $1\leq i \leq 3$ and
$4\leq j \leq 6$, and $R_i$ degenerates to $R'_i
+ \Pis_{jk}$ whenever $\{i,j,k\} = \{1,2,3\}$ or
$\{4,5,6\}$. Since we know (see~\cite{don}) that
each component of $F(X')$ is simple when $F(X')$
is viewed as a cycle on the grassmanian
$G(2,\PP^4)$, it follows that our $9$ planes and
$6$ rulings account for all of the components of
$F(X)$ (and each is simple on $F(X)$ when viewed
as a cycle on $G(2,\PP^4)$).
\end{proof}

For future use we summarize in the following
proposition the intersection pattern of the
components of $F(X)$:
\begin{proposition} {\bf 1.\ }
$\Pis_{14}\cap \Pis_{15} = \emptyset$.

\noindent {\bf 2.\ } $\Pis_{14}\cap \Pis_{25}$
consists of the one point corresponding to the
line $\Pi_{14}\cap \Pi_{25}$.

\noindent {\bf 3.\ } $\Pis_{14}\cap R_1$
consists of the one point corresponding to the
line $\ov{O_{15}O_{16}}$.

\noindent {\bf 4.\ } $\Pis_{14}\cap R_2$ is the
line of lines in $\Pi_{14}$ through $O_{34}$.

\noindent {\bf 5.\ } $R_1\cap R_2$ consists of
five points.

\noindent {\bf 6.\ } $R_1\cap R_4$ consists of
the conic $D$ of lines in $X$ passing through
$O_{36}$ and meeting $\Pi_{36}$, and the two
points corresponding to the lines
$\ov{O_{25}O_{14}}$ and $\ov{O_{15}O_{24}}$.
\end{proposition}

The other intersections are obtained by applying
an automorphism and a monodromy. We omit the
routine proof.

\section{The genus two case}
\nopagebreak

\subsection{More symmetry}
\label{more}

\nopagebreak When the base curve $C$ has genus
2, there is a group $(\ZZ/2\ZZ)^3$ extending the
symmetry group $(\ZZ/2\ZZ)^2$ which acts on
$\Cpm$ in other genera. To see this, we use the
following construction.

Start with three pairs of points
\[
a_{i,\epsilon}, \  i=1,2,3, \ \epsilon=0,1 \] in
$\PP^1$. Let $\PP^1_i$ be the double cover of
$\PP^1$ branched at the two points
$a_{i,\epsilon}, \ \epsilon=0,1$. The fiber
product $\PP^1_1 \times_{\PP^1} \PP^1_2
\times_{\PP^1} \PP^1_3$ is a curve ${_{5}}C$ of
genus 5. It admits a $(\ZZ/2\ZZ)^3$ action. The
seven level-1 quotients, or quotients by
subgroups $(\ZZ/2\ZZ)^2$, are the three
$\PP^1_i$ plus the three elliptic curves $E_i$
branched at the four points $a_{j,\epsilon}, \ j
\neq i, \ \epsilon=0,1$, and the genus 2
hyperelliptic curve ${_{2}}C$ branched at all
six points. The seven level-2 quotients, or
quotients by subgroups $\ZZ/2\ZZ$, consist of
the three elliptic curves $\widetilde{E_i} :=
\PP^1_j \times_{\PP^1} \PP^1_k$ (where $\{i,j,k
\} = \{1,2,3 \}$), plus the three genus 3 curves
${_{3}}C_i := \PP^1_i \times_{\PP^1} E_i$, and
one more genus 3 curve ${_{3}}C$ whose quotients
are the three $E_i$.

We can recover an instance of
diagram~\eqref{E:tower}) by relabelling:

$${_{5}}C = \Cpm, \ {_{3}}C_1=C_\qi, \
{_{3}}C_2=C_\qj, \ {_{3}}C_3=C_\qk, \
{_{2}}C=C,$$ and choosing a double cover
$\widetilde{C} \to \Cpm$ such that
$\widetilde{C}$ is Galois over each ${_{3}}C_i$
with group $\ZZ/4\ZZ$. But in fact, any
diagram~\eqref{E:tower} with $g=2$ arises this
way, and uniquely. The point is that the base
curve ${_{2}}C$ is hyperelliptic, and the
hyperelliptic involution acts on $J({_{2}}C)$ as
$-1$, so it preserves all points of order 2 and
all double covers. In fact, any double cover
such as $C_\qi \to {_{2}}C$ is Galois over
$\PP^1$ with group $(\ZZ/2\ZZ)^2$ and quotients
${_{2}}C, E_\qi, \PP^1_\qi$ of genera $2,1,0$
respectively. In particular, this gives three
double covers $\PP^1_\qi$, $\PP^1_\qj$,
$\PP^1_\qk$ of $\PP^1$. If we relabel them
$\PP^1_i, \ i=1,2,3, $ we are back in the
situation of the previous paragraph.

\subsection{Even more symmetry}
\label{evenmore}

Now start with a set $S$ of five points in
$\PP^1$. We label the five distinct elements, in
any order, as $i,j,k,l,m.$ There are
$15=2^{(5-1)}- 1$ non-empty even subsets of $S$,
giving $15$ branched double covers of $\PP^1$.
These are the $15$ level-1 quotients of their
common Galois closure, a curve ${_{5}}C$ which
is Galois over $\PP^1$ with group
$(\ZZ/2\ZZ)^4$. We enumerate all the quotients
of ${_{5}}C$:

\begin{tabular}{|ccc|l|}
\multicolumn{3}{l}{\bf Level 1}\\ \hline 10 &
rational & curves & $\PP^1_{i,j}$ is branched at
2 points $i,j \in S$ \\ \hline 5 & elliptic &
curves & $E_i$ is branched at 4 points $S
\setminus i$
\\ \hline\hline
\multicolumn{3}{l}{\bf Level 2} \\ \hline 10 &
rational & curves & $\PP^1_{i,j,k}$ has
quotients $\PP^1_{i,j}, \PP^1_{i,k},
\PP^1_{j,k}$\\ \hline 15 & elliptic & curves &
$E_{ij,kl} := \PP^1_{i,j} \times _\PP^1
\PP^1_{k,l}$ \\ \hline 10 & genus 2 & curves &
${_{2}}C_{i,j}$ has quotients $\PP^1_{i,j},
E_i,E_j$ \\ \hline\hline \multicolumn{3}{l}{\bf
Level 3} \\ \hline 5 & elliptic & curves &
$\widetilde{E}_i$
has level-1 quotients: \\
 & & & $E_i, \PP^1_{j,k}, \PP^1_{j,l},
 \PP^1_{j,m}, \PP^1_{k,l},
\PP^1_{k,m}, \PP^1_{l,m}, $   \\
& & & and level-2 quotients: \\
 & & & $\PP^1_{j,k,l}, \PP^1_{j,k,l},
 \PP^1_{j,k,l}, \PP^1_{j,k,l},
E_{jk,lm}, E_{jl,km}, E_{jm,kl}$ \\ \hline 10 &
genus 3 & curves & ${_{3}}C_{i,j}$
has level-1 quotients: \\
 & & & $\PP^1_{i,j}, \PP^1_{k,l},
 \PP^1_{k,m}, \PP^1_{l,m}, E_k, E_l,
E_m$  \\
& & & and level-2 quotients: \\
 & & & $\PP^1_{k,l,m}, E_{ij,kl},
 E_{ij,km}, E_{ij,lm},
{_{2}}C_{k,l}, {_{2}}C_{k,m}, {_{2}}C_{l,m}.$
\\ \hline\hline
\multicolumn{3}{l}{\bf Level 4} \\ \hline 1 &
genus 5 & curve & ${_{5}}C$ \\ \hline \hline
\end{tabular}
\vspace{0.3cm}

Note that $\widetilde{E}_i$ is Galois over $E_i$
with group $(\ZZ/2\ZZ)^2$ and intermediate
covers $E_{jk,lm}$, $E_{jl,km}$, and
$E_{jm,kl}$. These are all the double covers of
$E_i$; so $\widetilde{E}_i$ is isomorphic to
$E_i$, and the degree 4 map $\widetilde{E}_i \to
E_i$ is multiplication by 2.

The Galois group of ${_{5}}C$ over $\PP^1_{ij}$
is $(\ZZ/2\ZZ)^3$. In $\PP^1_{ij}$ we have six
branch points in three pairs, namely the inverse
images of $S \setminus \{i,j\}$.  So the curve
${_{5}}C$ can be viewed as our previous $\Cpm$
in ten distinct ways, over the ten rational
curves $\PP^1_{ij}$ and the corresponding
genus-2 base curves ${_{2}}C_{i,j}$.

In this special case we can also describe
quaternion covers $\widetilde{C} \to \Cpm$ quite
explicitly. Let $q$ be either of the two points
of $E_m$ above $m \in S \subset \PP^1$. Its
inverse image in $\widetilde{E}_m$ is the set of
four points $q_a, \ a=1,2,3,4$ satisfying
$2q_a=q$. Now on $\widetilde{E}_m$ we have a
natural line bundle $\cL_m \in
{\Pic}^2(\widetilde{E}_m)$ such that $\cL_m
^{\otimes 2}$  has a section $s$ vanishing at
the four points $q_a$. Namely, $\cL_m$ is
isomorphic to $\cO_{\widetilde{E}_m}(2q_a)$, for
any $a$. The inverse image in $\cL_m$ of the
section $s$, under the squaring map, gives a
double cover ${_{3}}C_{m} \to \widetilde{E}_m$
branched at the four points $q_a$. Explicitly,
if we write the equation of $\widetilde{E}_m$ as
a double cover of $\PP^1_{i,j,k}$ as:
$$y^2 = \Pi_{a=1}^4 (x-\lambda_a),$$
with $q_a$ the point with coordinates
$(x=\lambda_a, y=0),$ then ${_{3}}C_{m}$ has
equation
$$y^4 = \Pi_{a=1}^4 (x-\lambda_a).$$
In particular, ${_{3}}C_{m}$ is
$\ZZ/4\ZZ$-Galois over $\PP^1_{i,j,k}$. It
follows that the fiber product:
$$\widetilde{C} := {_{5}}C
\times _ {\widetilde{E}_m} {_{3}}C_m =
{_{3}}C_{l,m}  \times _
{\PP^1_{ijk}}{{_{3}}C_m}$$ is a
$\ZZ/4\ZZ$-Galois cover of ${_{3}}C_{l,m}$.
Similarly, this same $\widetilde{C}$ is also a
$\ZZ/4\ZZ$-Galois cover of ${_{3}}C_{i,m}$,
${_{3}}C_{j,m}$, and ${_{3}}C_{k,m}$. In
particular, $\widetilde{C}$ is quaternionic over
${\PP^1_{ijk}}$ (and also over ${\PP^1_{ijl}},
{\PP^1_{ikl}}$, and ${\PP^1_{jkl}}$).

We note that what we get this way is a special
case of the general genus 2 quaternionic
towers~\eqref{E:tower}: the general case depends
on 3 parameters, while this special case depends
on only two parameters. The curves ${_{5}}C$ in
this two dimensional family are known as Humbert
curves, cf. \cite{don,var}. Varley shows \cite
{var} that the covers $\widetilde{C} \to
{_{5}}C$ all have the same Prym, a certain
4-dimensional non-hyperelliptic ppav with 10
vanishing theta nulls.

\subsection{Abelian fourfolds and cubic
threefolds} \label{afct} We need to recall some
features of the Prym map in genus $5$. Our
references in this subsection are \cite{don1,
don}. Let $\cA_g$ be the moduli space of
g-dimensional ppav's, and $\cR\cA_g$ the moduli
space of g-dimensional ppav's with a marked
point of order 2. Let $\cM_g$ be the moduli
space of curves of genus g, and $\cR_g$ the
moduli space of curves with a marked point of
order 2 in their Jacobian. Let $\cC$ be the
moduli space of cubic threefolds whose only
singularities are some ordinary double points.
There is a corresponding moduli space $\cR\cC$
of cubic threefolds together with a point of
order 2 in their intermediate Jacobian. In fact,
this space splits into even and odd components:
$\cR\cC = \cR\cC^+ \cup \cR\cC^-$, distinguished
by an appropriate $\Ztwo$-valued function.
Similarly, let $\cQ$ be the moduli space of
plane quintic curves $Q$ whose only
singularities are some ordinary double points.
There is a corresponding moduli space $\cR\cQ$
of plane quintic curves together with a point of
order 2 in their compactified Jacobian. Again,
this space splits into even and odd components:
$\cR\cQ = \cR\cQ^+ \cup \cR\cQ^-$, distinguished
by an appropriate $\Ztwo$-valued function.

One of the basic results about the Prym map:

$$\cP: \cR_5 \to \cA_4$$

is that it factors through a rational map:

$$\kappa: \cR_5 \to \cR\cC^+$$

followed by a birational isomorphism:

$$\chi: \cR\cC^+ \to \cA_4.$$

These are constructed as follows. For more
details, see~\cite{don}.

\begin{itemize}

\item A pair $(X,l)$, where $X \subset \PP^4$
is a cubic threefold (with nodes at worst) and
$l \subset \PP^4$ is a line contained in $X$,
determines a plane quintic curve $Q$ and an
(odd) double cover $\widetilde{Q}_{\sigma} \to
Q$. Explicitly, $\widetilde{Q}_{\sigma}$
parameterizes the lines in $X$ meeting $l$,
while $Q$ parameterizes the planes through $l$
which intersect $X$ residually in two additional
lines.

\item Conversely, given $Q$ and an odd $\sigma$,
we can recover the pair $(X,l)$. When $Q$ is
non-singular, $X$ is characterized as the unique
cubic threefold whose intermediate Jacobian is
isomorphic to $\Prym (Q,\sigma)$. In case $Q$ is
singular, we describe an explicit construction
below, in the proof of the implication {\bf
$(2)\Rightarrow(1)$} in Theorem \ref{detailed}.

\item A non-hyperelliptic curve $C \in \cM_5$
determines a plane quintic curve $Q$ and an
(even) double cover $\widetilde{Q}_{\nu} \to Q$
such that $\cP(Q, \nu) \cong \Jac(C).$
Explicitly, $\widetilde{Q}_{\nu}$ is the
singular locus of the theta divisor of
$\Jac(C)$, so it parameterizes linear systems
$g^1_4$ on $C$, while $Q$ is the quotient of
$\widetilde{Q}_{\nu}$ by the involution $-1$ of
$\Jac(C)$, and it parameterizes quadrics of rank
4 in $\PP^4$ through the canonical image of $C$.

\item In the above situation, a point of order
2: $\mu \in \Jac(C)[2]$, determines via
Mumford's isomorphism (cf. \cite{mum2} or
\cite[Theorem 1.4.2]{don}), a pair of points of
order two: $\sigma, \nu\sigma \in \Jac(Q)[2]$.
One of these, say $\nu\sigma$, is even, while
the other, $\sigma$, is odd.

\item So, given $(C,\mu) \in \cR_5$, we set
$\kappa(C,\mu)=(X,\delta)$, where $X$ is the
cubic threefold corresponding to $(Q,\sigma)$,
and $\delta$ is the image of $\nu$ under
Mumford's isomorphism. It is automatically even.

\item Finally, given $(X,\delta) \in \cR\cC^+$,
choose a line $l$ in $X$. The pair $(X,l)$
determines the quintic $Q$ and its odd double
cover $\widetilde{Q}_{\sigma} \to Q$, while
$\delta$ determines a second cover
$\widetilde{Q}_{\nu} \to Q$. Then $\Prym(Q,\nu)$
is the Jacobian of a curve $C$ which can be
described explicitly, and $\sigma$ descends to a
point $\mu \in \Jac(C)[2]$. We then set
$\chi(X,\delta) := \Prym(C, \mu)$. The result
turns out to be independent of the choice of
$l$, by the tetragonal construction.
\end{itemize}

\subsection{The main results}

Our main result is that under the correspondence
$\chi$, the four dimensional quaternionic
abelian varieties correspond to the nine-nodal
cubic threefolds. Each of these nine-nodal cubic
threefolds comes equipped with a natural,
"allowable" point of order two. Before giving
the precise statement of the theorem, we need to
explain this lift.

Let $X$ be a nine-nodal cubic threefold. We use
the notation of section \ref{cubics}. Let
$R=R_6$ be a ruling consisting of all lines in
$X$ meeting three transversal planes $\Pi_{16},
\Pi_{26}, \Pi_{36}$, and let $l \in R$ be one
such line. This determines a plane quintic
$Q=Q_l$ and its double cover
$\widetilde{Q}_{\sigma}$. We note that
$\widetilde{Q}_{\sigma}$ contains the three
pencils $L_a^0$ of lines in the plane $\Pi_{a6}$
which pass through the point $\Pi_{a6} \cap l$,
for $a \in A = \{1,2,3\}$. Therefore the plane
quintic $Q$ contains three lines $L_a$, and
residually a conic $D$. It follows that the
cover $\widetilde{Q}_{\sigma} \to Q$ is \'etale:
$\widetilde{Q}_{\sigma} = (\cup_{a=1}^3
\cup_{\epsilon=0}^1 L_a^{\epsilon}) \cup D^0
\cup D^1,$ where each $D^{\delta}$ meets each
$L_a^{\epsilon}$   in one point.

Our $X$ represents a point of the moduli space
$\cC$. A lift of $X$ to $\cR \cC^+$ is
determined by a second double cover
$\widetilde{Q}_{\nu}$ of $Q$, with class $\nu$
which is orthogonal to $\sigma$ and even. There
is a natural choice for such a double cover, and
hence for the lift, namely the unique allowable
one: $\widetilde{Q}_{\nu}$ is the unique cover
which is branched over all nine nodes of $Q$.
Explicitly, $\widetilde{Q}_{\nu} = E_1 \cup E_2
\cup E_3 \cup {{_{2}}\widetilde{D}}$, where $E_a
\to L_a$ is a double cover branched at the four
points where $L_a$ meets the other components,
and ${{_{2}}\widetilde{D}} \to D$ is a double
cover branched at the six intersection points of
the conic with the lines.

\begin{lemma}\label{parity}
The allowable cover $\widetilde{Q}_{\nu} \to Q$
constructed from a pair $(X,l)$ is even and
orthogonal to $\sigma$.
\end{lemma}

We can  now state our main results.

\begin{theorem}\label{main}
The correspondence $\chi$ takes the nine-nodal
cubic threefolds (with their unique allowable
lift to $\cR \cC^+$) to the four dimensional
quaternionic abelian varieties.
\end{theorem}

This follows immediately from the following more
detailed version:

\begin{theorem}\label{detailed}
The following data are equivalent:
\begin{enumerate}

\item Pairs $(X,l)$ where $X$ is a nine-nodal
cubic threefold and $l$ is a line in a ruling
$R$ on $X$.

\item Pairs $(Q,\widetilde{Q}_{\sigma})$
where $Q=L_1 \cup L_2 \cup L_3  \cup \Delta$ is
a reducible quintic consisting of three lines
and a conic, and $\widetilde{Q}_{\sigma} \to Q$
is an \'etale double cover,
$\widetilde{Q}_{\sigma} = (\cup_{a=1}^3
\cup_{\epsilon=0}^1 L_a^{\epsilon}) \cup
\Delta^0 \cup \Delta^1,$ where each
$\Delta^{\delta}$ meets each $L_a^{\epsilon}$ in
one point, and $L_a^{\epsilon}$ meets
$L_{a'}^{\epsilon'}$ if and only if $a \neq a'$
and $\epsilon \neq \epsilon'$.

\item A curve $\widetilde{C} \in \cM_9$ with a
fixed-point free action of the quaternion group
$\quat$.

\end{enumerate}
\end{theorem}

\begin{proof} (of Theorem \ref{detailed} and
Lemma \ref{parity})

\vskip 0.5cm \noindent {\bf $(1)\Ra(2)$:}

We saw above how to go from $(X,l)$ to a
reducible plane quintic $Q=L_1 \cup L_2 \cup L_3
\cup \Delta$ and an \'etale double cover
$\widetilde{Q}_{\sigma} = (\cup_{a=1}^3
\cup_{\epsilon=0}^1 L^{\epsilon}_a) \cup
\Delta^0 \cup \Delta^1.$ The intersection
properties of the components of
$\widetilde{Q}_{\sigma}$ can be determined
directly from the explicit
formula~\eqref{E:nine}. An alternative is to
return to the degeneration used in Lemma
\ref{Fanos}, in which $X$ goes to the Segre
cubic and $l\in R_6$ goes to $l'\in R'_6$. The
cover $\widetilde{Q'}_{\sigma} \to Q'$
corresponding to $(X',l')$ is easy to determine,
because of the larger symmetry present in this
case. It was described, for example, in
\cite{don}, formula (5.17.4):

$$Q' = (\cup_{i=1}^5 L_i^{\epsilon}) /
(p_{i,j} \sim p_{j,i}, \ i \neq j),$$

\begin{equation}
\label{aaa} \widetilde{Q'}_{\sigma} =
(\cup_{i=1}^5 \cup_{\epsilon=0}^1
L_i^{\epsilon}) / (p_{i,j}^0 \sim p_{j,i}^1, \ i
\neq j).
\end{equation}

Under our degeneration, the conic $\Delta$
splits into $L_4 \cup L_5$. The cover
$\widetilde{Q}_{\sigma} \to Q$ given in the
theorem is the only one which specializes
correctly.

\vskip 0.5cm \noindent {\bf $(2)\Ra(1)$:}

To go in the opposite direction, consider first
the more general situation, where we start with
a pair $(Q,\widetilde{Q}_{\sigma})$, where $Q$
is any quintic with at least one node $o$ over
which $\widetilde{Q}$ is \'etale. We can
explicitly exhibit the corresponding cubic
threefold $X$ and line $l$ as follows.
Projection from $o$ shows that the partial
normalization $T$ of $Q$ at $o$ is a trigonal
curve of arithmetic genus 5, with a double cover
$\widetilde{T}$ obtained by normalizing
$\widetilde{Q}$ above $o$. The trigonal
construction takes the pair $T,\widetilde{T}$ to
a curve $B$ of genus 4 which comes equipped with
a $g^1_4$ linear system. The canonical map sends
$B$ to $\PP^3$, and the homogeneous ideal of the
image is generated by a quadric $f_2$ and a
cubic $f_3$. The inhomogeneous equation $f_2
+f_3=0$ then determines a Zariski open piece of
our cubic $X$ as a hypersurface in affine
4-space, and $X$ is recovered as the closure in
$\PP^4$. Hence $B$ can be naturally identified
with the curve $C_{O,X}$ of lines on $X$ through
$O$ (introduced in Theorem~\ref{nine}). By
Proposition~\ref{Fanos}, the Fano surface $F(X)$
parameterizing lines in $X$ can be described as
the symmetric product $S^2B$ modulo certain
identifications. The $g^1_4$ linear system on
$B$ is necessarily of the form $\omega_B(-p-q)$
for two points $p,q$ in (the smooth part of)
$B$, where $\omega_B$ is the canonical bundle.
We then recover $l$ as the line corresponding to
the point of $F(X)$ given by the image of $p+q
\in S^2B$. From section (5.11.2) of \cite{don}
it follows that this construction is indeed
inverse to our construction of
$(Q,\widetilde{Q}_{\sigma})$ from $(X,l)$.

Returning to our special case, we now see that
it merely remains to check that this line $l$
lies indeed on a ruling (and not on a dual
plane). Assume (as we may) that, in the previous
notation, $O=O_{3,6}$. Then we will show a more
precise result:
\begin{claim}
The line $l$ is on the ruling $R_3$ or $R_6$ if
and only if the node $o$ is on the intersection
of two lines $L_i$, $L_j$ of $Q$; on the other
hand, $l$ is on one of the rulings $R_1$, $R_2$,
$R_4$ or $R_5$ if and only if the node $o$ is on
the intersection of a line of $Q$ and the conic
$D$. (The monodromy action permutes all the
cases of a given type.)
\end{claim}
To see this we use the analysis
in~\eqref{E:sym2}: this tells us which
components of $B$ must contain the points $p$,
$q$ in order for the line $l$ to be on a given
ruling $R_i$:
\begin{equation}
\label{pointcomp}
\begin{array}{ll}
l\in R_3 &\Leftrightarrow
p\in L_{34},\;\;q\in L_{35}\\
l\in R_6 &\Leftrightarrow
p\in L_{16},\;\;q\in L_{26}\\
l\in R_i \quad\mbox{for}\quad i=1,2
&\Leftrightarrow p\in L_{i6},\;\;q\in D\\
l\in R_j \quad\mbox{for}\quad j=4,5
&\Leftrightarrow p\in L_{3j},\;\;q\in D.
\end{array}
\end{equation}
As in the general case, the degree $4$ map
$\pi:B\ra \PP^1$ is given by projecting from the
line $\overline{pq}\subset \PP^3$. Here this
projection has degree $0$ (i.e. it is constant)
on the line components of $B$ through $p$, $q$;
on the remaining line components of $B$ the
degree of $\pi$ is $1$, and on $D$ it is $1$
when $q$ lies on $D$ (``the second case'') and
$2$ otherwise (``the first case''). To prove the
claim, we must show that the degree $4$ map
$\pi:B\ra \PP^1$ arises from the double cover
$\widetilde{Q}_{\sigma} \to Q$ of the trigonal
curve $\pi':Q\ra \PP^1$ by the trigonal
construction. There are two cases to examine,
namely when $\pi'$ is projection from the
intersection of two lines (``the first case''),
and when it is projection from a point of
intersection of $\Delta$ and a line ``the second
case''). As before, $\pi'$ has degree $0$ on the
lines through the center of projection, and it
is straightforward to see from the definition of
the trigonal construction, that the two cases we
distinguished for $\pi'$ yield the respective
cases we distinguished for $\pi$.

\vskip 0.5cm \noindent {\bf $(3)\Ra(2)$:} Recall
first that by part (3) of Lemma~\ref{liftpsi}
there are $4$ possible double covers $\Ctil$\/'s
covering a given $\Cpm$ in (3). Similarly, there
are $4$ possible double covers
$\widetilde{Q}_\sigma$ in (2): for each
$a=1,2,3$ we must choose which of the two points
of $\Delta^0$ which lie above the points where
$\Delta$ and $L_a$ intersect is on $L_a^0$. Of
the resulting $2^3=8$ possibilities each choice
is isomorphic with the ``opposite'' one,
obtained by interchanging $\Delta^0$ with
$\Delta^1$ and each chosen point of $\Delta\cap
L_a$ with the other one. We now claim that the
parameter spaces, $\cR_3$ for the coverings
$\Ctil/\Cpm$ and $\cR_2$ for the
$\widetilde{Q}_\sigma/Q$\/'s, form irreducible
spaces. For $\cR_3$ this is
Proposition~\ref{strata}. For $\cR_2$ notice
first that the space of $Q$\/'s (conics and
three lines) is manifestly irreducible. The same
is then true for the allowable covers
$\widetilde{Q}_\nu$ of $Q$ by their uniqueness.
Moreover, monodromy allows us to ``turn around''
individually each of the lines so its points of
intersection with the conic are interchanged.
This shows that the $4$ covers
$\widetilde{Q}_\sigma$\/'s of a given $Q$ are in
the same component, proving the claim.

Now suppose that we are given the curve $\Ctil$
with an action of $G$, hence the quotient $\Cpm$
and the entire tower~\eqref{E:tower}. Let
$\widetilde{Q}_{\nu} \to Q$ be the quintic
double cover corresponding to the genus-5 curve
$\Cpm$, and let $\widetilde{Q}'_{\sigma} \to Q$
be the double cover inducing $\widetilde{C} \to
C_\pm$ via Mumford's isomorphism,
cf.~\cite{mum2} or \cite[Theorem 1.4.2]{don}. As
we saw in subsection \ref{more}, $\Cpm$ has
three elliptic quotient curves $\widetilde{E_i},
\ i=1,2,3$. It follows that
$\widetilde{Q}_{\nu}$, which parameterizes
linear systems $g^1_4$ on $\Cpm$, contains three
elliptic curves, which can be canonically
identified with the Picard varieties
$\Pic^2(\widetilde{E_i})$. Therefore, $Q$
contains three lines $L_i=\PP^1_i$, and
residually a conic $\Delta$. We claim that the
double cover $\widetilde{Q}'_{\sigma} \to Q$ is
one of the double covers $\widetilde{Q}_{\sigma}
\to Q$ described in part (2) of the theorem.

This is known to be true after we specialize the
general curves $\Cpm$ of subsection \ref{more}
to the Humbert curves of \ref{evenmore}: the
conic $\Delta$ breaks further to two lines
$L_4,L_5$, so that each of $\Delta^0$ and
$\Delta^1$ breaks into two lines $L^0_4\cup
L^1_5$ and $L^1_4\cup L^0_5$ respectively. The
double cover obtained from the Segre cubic is
specified in (\ref{aaa}) and agrees with the
double cover $\widetilde{Q}_{\sigma} \to Q$
described in part (2) of the theorem.

We now claim that the same must hold in general,
namely that for every $\Ctil$ in (3) the
covering $\Ctil_\sigma$ we obtained is one of
the $4$ covers given in (2). Indeed, let $\cQ$
be the irreducible variety parameterizing the
reducible quintics as in part (2) of the
theorem, and let $\cR\cQ \to \cQ$ be the \`etale
cover parameterizing all \`etale double covers
of such quintics. We are given two irreducible
subcovers $\cR_2 \to \cQ$ and $\cR_3 \to \cQ$ of
$\cR\cQ \to \cQ$, parametrizing the \`etale
double covers coming from (2) and (3)
respectively. (We noted that each is a
four-sheeted cover.) Now two irreducible
subcovers of an \`etale cover must either be
disjoint or coincide. But our $\cR_2$ and
$\cR_3$ intersect at the Segre point
$\widetilde{Q}_{\sigma,0}/Q_0$, at which $Q=Q_0$
consists of five lines. They therefore coincide
as claimed.

\vskip 0.5cm \noindent {\bf $(2)\Ra(3)$:}

As we saw above, the allowable double cover
$\widetilde{Q}_{\nu} \to Q$ is uniquely
determined by $Q$. We recover $\Cpm$ as the
unique curve whose Jacobian is isomorphic (as a
ppav) to $\Prym(\widetilde{Q}_{\nu} / Q)$, and
$\widetilde{C} \to \Cpm$ is the double cover
induced via Mumford's isomorphism from
$\widetilde{Q}'_{\sigma} \to Q$. This is clearly
the inverse of the previous construction, so we
are done with the theorem.

If we start with data (3), the construction   of
$\widetilde{Q}_{\sigma} \to Q$ involves
Mumford's isomorphism, so the cover
$\widetilde{Q}_{\sigma} \to Q$ is automatically
orthogonal to $\widetilde{Q}_{\nu} \to Q$.
Moreover, this  cover $\widetilde{Q}_{\nu} \to
Q$ is even, since its Prym is the Jacobian of a
curve $\Cpm$ (rather than the intermediate
Jacobian of a cubic threefold). By the theorem,
it follows that the same holds if we start with
data (1), proving the lemma.

\end{proof}

Now that we have a fairly complete description
of the Prym map on the parameter space
$\cM(2;0,0,0)$ of proposition \ref{strata}, it
is easy to specialize further and find its
behavior on the boundary strata $\cM(1;2,0,0)$
and $\cM(0;2,2,0)$. We work out the latter in
some detail.

We want the genus 2 curve $C$ at the bottom of
tower~\eqref{E:tower} to degenerate, within a
fixed fiber of the Prym map, to a rational curve
with two nodes, and we want to trace what
happens to $\Cpm, \widetilde{C}$ in this limit.
This can be easily arranged, in the language of
subsection \ref{more}, for example by letting
the branch points $a_{i,0}, a_{i,1}$ coincide
for $i=2,3$. We relabel the surviving points
$a_{1,0}, a_{1,1}, a_2, a_3$. At level 1 of the
$(\Ztwo)^3$-diagram we find that $\PP^1_1$
remains a smooth rational curve, doubly covering
$\PP^1$ with  branch points $a_{1,0}, a_{1,1}$;
but $\PP^1_i$ for $i=2,3$ degenerates to a
reducible curve, consisting of two copies
$\PP^1_{i,\epsilon_i}$ of $\PP^1$, $\epsilon_i
\in (\Ztwo)$, intersecting each other above
$a_i$. It follows immediately that ${_5}C =
\Cpm$ has four components $C_{\epsilon_2,
\epsilon_3}, \ \epsilon_2, \epsilon_3 \in
\Ztwo$. Each component $C_{\epsilon_2,
\epsilon_3}$ is isomorphic, as a double cover of
$\PP^1$, to $\PP^1_1$, i.e. it is branched at
$a_{1,0}, a_{1,1}$. Component $C_{\epsilon_2,
\epsilon_3}$ meets component $C_{1+\epsilon_2,
\epsilon_3}$ in the two points above $a_2$; it
meets component $C_{\epsilon_2, 1+\epsilon_3}$
in the two points above $a_3$; and it does not
meet component $C_{1+\epsilon_2, 1+\epsilon_3}$.
Finally, up to isomorphism there is only one
quaternionic cover $\widetilde{C} \to \Cpm$,
namely the unique allowable cover of $\Cpm$. Let
$E(a_{1,0}, a_{1,1}, a_2, a_3)$ be the elliptic
curve which is the double cover of $\PP^1_1$
branched at the four points above $a_2, a_3 \in
\PP^1$. Then $\widetilde{C}$ consists of four
copies of the same elliptic curve $E(a_{1,0},
a_{1,1}, a_2, a_3)$ glued at their ramification
points. The Prym is then isogenous to
$(E(a_{1,0}, a_{1,1}, a_2, a_3))^{\times 4}$.

In a special case, this degeneration picture was
obtained in Remark 3 on the last page of
\cite{var}; in fact, Varley's situation is
precisely the case of our subsection
\ref{evenmore}, where two pairs among the five
points $S \subset \PP^1$ coalesce, say $i$ and
$l$ go to $0$ while $j$ and $m$ go to $\infty$
and $k$ is at $1$. The six points on
$\PP^1_{i,j}$ then  coincide as in the previous
paragraph: $a_{1,0}=1, a_{1,1}=-1, a_2=0,
a_3=\infty$, showing that the Prym of an Humbert
curve is isogenous to the fourth power of the
harmonic elliptic curve. In our more general
setting, the elliptic curve $E(a_{1,0}, a_{1,1},
a_2, a_3)$ is arbitrary: we can take for
instance $a_{1,0}=0, a_{1,1}=\infty$. We then
see that $E(0, \infty, a_2, a_3)$ is the double
cover of $\PP^1_1$ branched at $\pm \sqrt{a_2},
\pm \sqrt{a_3}$, and this has  variable modulus.
We conclude:

\begin{corollary}
\label{end} The Prym map $\Phi$ sends
$\cM(2;0,0,0)$ and each of the spaces
$\cM(1;2,0,0)$ and $\cM(0;2,2,0)$ onto the
Shimura curve $\Shim$ parameterizing
4-dimensional ppav's $A$ with quaternionic
multiplication. Each 4-dimensional ppav $A$ with
quaternionic multiplication is isogenous to the
fourth power $E^{\times 4}$ of some elliptic
curve $E$, and every $E$ occurs for some $A$.
(The precise isogeny is given in
Corollary~\ref{PEL2}.)
\end{corollary}

\section{Appendix}

In this appendix we will sketch the proof of the
following result, mentioned in the Introduction.
\begin{theorem}
Let $T^{\ast\ast}/T$ be a cyclic, $4$\/-sheeted
Galois unramified cover of a (smooth projective
complex irreducible) general trigonal curve $T$
of genus $g>1$. Identify $\Gal(T^{\ast\ast}/T)$
with $\left<\qi\right>$, and let $T^\ast/T$ be
the intermediate $2$\/-sheeted cover. Then the
$\left<\qi\right>$\/-action on $P =
\Prym(T^{\ast\ast}/T^\ast)$ does not extend to a
$G$\/-action.
\end{theorem}
\begin{proof}
The locus of hyperelliptic curves in $\cM_g$ is
in the closure of the trigonal locus. Hence
existence of a $G$\/-action in the trigonal case
implies the same for each hyperelliptic curve
$H$. In the trigonal case there is only one type
of cyclic covers (under monodromy), but over the
hyperelliptic locus there are several types of
unramified double covers $H^\ast/H$ and
a-fortiori of $4$\/-sheeted cyclic covers
$H^{\ast\ast}/H$. A double cover $H^\ast/H$ is
determined by the choice of a subset $D$ of even
cardinality of the Weierstrass points of $H$, up
to replacing this subset by its complement. The
irreducibility of the space of $T^{\ast\ast}/T$
over the trigonal locus implies that the
$G$\/-action must then exist for each such
(hyperelliptic) type. Thus it suffices to
consider the easiest type, when the set $D$
consists of two Weierstrass points, which we
take as the points of $H$ over $0,\infty$ in
$\PP^1$. In terms of a coordinate $x$ on
$\PP^1$, the curve $H$ is given by an affine
equation $y^2=xf(x)$, where $f$ has degree $2g$
and simple roots. Then $H^\ast$ is given by
$(y/u)^2=f(u^2)$ (with $x=u^2$). In particular,
$H^\ast$ is hyperelliptic. We will need the
following two lemmas:

\begin{lemma}
\label{l1} In this case $P$ is the square of a
general hyperelliptic jacobian $J(C)$ (up to
isogeny).
\end{lemma}
\begin{lemma}
\label{l2} For a general hyperelliptic curve $C$
we have $\End(\Jac(C)) = \ZZ$.
\end{lemma}
Assuming Lemmas~\ref{l1},~\ref{l2} it is easy to
conclude the proof of the Theorem. Indeed, the
$\QQ$\/-endomorphism ring of $P$ is then
$\Mat_{2\times 2}(\QQ)$. The group $G$ must then
embed into $\GL_2(\QQ)$ (the action of $G$ on
$P$ must be faithful since $-1\in G$ acts as
$-1$ on $P$). This is a contradiction, because
$G$ does not embed even into $\GL_2(\RR)$.
\end{proof}

In the proof of Lemma~\ref{l1} we will need the
following third Lemma:
\begin{lemma}
\label{cyclic}
 Let $F$ be a curve, let $F^\ast$
be an unramified double cover of $F$
corresponding to the $\sigma\in\Jac F[2]$, and
let $F{\ast\ast}$ be an unramified double cover
of $F^\ast$ corresponding to the class
$\sigma^\ast\in\Jac F^\ast[2]$. Then
$F^{\ast\ast}$ is cyclic Galois over $F$ if and
only if the norm map $\Nm:\Jac F^\ast \ra \Jac
F$ maps $\sigma^\ast$ to $\sigma$.
\end{lemma}
\begin{proof} We will prove this presumably
well-known Lemma~\ref{cyclic} since we do not
know a reference for it. Assume first that
$F^{\ast\ast}/F$ is cyclic (of order $4$\/). As
in the proof of Proposition~\ref{strata} we
present the fundamental group
$\pi=\pi_1(F,\dot)$ of $F$ as generated by a
standard (symplectic) basis
$\alpha_1,\beta_1,\dots,\alpha_g,\beta_g$. We
assume, as we may, that $F^{\ast\ast}$ is the
cover corresponding to the kernel of the
homomorphism $\pi\ra\ZZ/4\ZZ$ given by sending
$\alpha_1$ to $1$ (in $\ZZ/4\ZZ$\/) and the
other generators to $0$. Then we may think of
$F$ as a genus one curve $E$, with fundamental
group generated by $\alpha_1$ and $\beta_1$, to
which a ``tail'' $T$ of genus $g-1$ is attached.
In this model $F^{\ast\ast}$ can be viewed as a
genus one curve $E^{\ast\ast}$, with fundamental
group generated by $\alpha^4\alpha_1$ and
$\beta_1$, with four copies of $T$ attached. The
deck transformations of $F^{\ast\ast}/F$ rotate
$E$ cyclically by quarter-turns along
$4\alpha_1$ and permute cyclically the four
$T$\/-tails. The intermediate cover $F^\ast$ is
then similarly a genus one curve $E^\ast$, whose
fundamental group is generated by $\alpha^\ast =
2\alpha_1$ and $\beta_1$, to which two copies of
$T$ are attached, with the deck transformation
rotating by half turns along $2\alpha_1$ while
permuting these two copies of $T$. This is
summarized by the following diagram:
\begin{equation}
\label{cyclicc}
\begin{array}{rlcccl}
\circlearrowleft\phantom{-} &T&&\curvearrowleft&&\\
&\,|&&&&\\
T-&E^{\ast\ast}-T&\qquad\ra\qquad&
T-E^\ast-T&\qquad\ra\qquad&E-T\\
&\,|&&&&\\ &T&&&&
\end{array}
\end{equation}

Under the duality between $\Jac(F)[2]$ and
$H_1(F,\ZZ/2\ZZ)$ we readily see that $\sigma$
is dual to $\alpha_1$. Indeed
$H_1(F,\ZZ)=H_1(E,\ZZ)\oplus H_1(T,\ZZ)$, with
$H_1(E,\ZZ)=\ZZ\alpha_1\oplus\ZZ\beta_1$.
Likewise  $H_1(F^\ast,\ZZ)=H_1(E^\ast,\ZZ)\oplus
H_1(T,\ZZ)^2$, with
$H_1(E^\ast,\ZZ)=\ZZ\alpha^\ast\oplus\ZZ\beta_1$.
The projection from $F^\ast$ to $F$ visibly maps
$H_1(F^\ast,\ZZ)$ onto the subgroup
$2\ZZ\alpha_1\oplus\ZZ\beta_1\oplus H_1(T,\ZZ)$
of $H_1(F^\ast,\ZZ)$. Reducing the coefficients
modulo $2$ we see that the annihilator of this
image modulo $2$ in $H_1(F,\ZZ/2\ZZ)$, under the
intersection pairing, is indeed (the class
modulo $2$ of) $\alpha_1$, which is therefore
dual to $\sigma$. Likewise $\alpha^\ast$ is dual
to $\sigma^\ast$. But the norm map above is dual
to the pull-back from $H_1(F,\ZZ/2\ZZ)$ to
$H_1(H^\ast,\ZZ/2\ZZ)$, which indeed maps the
class of $\alpha_1$ to the class of
$\alpha^\ast$. This proves the``only if'' part
of our claim.

To prove the ``if'' part (which is our main
concern), observe that the space $\cR_g^{(n)}$
of cyclic unramified $n$\/-sheeted covers of
curves of genus $g$ is connected by the same
argument as in the proof of
Proposition~\ref{strata}. It forms an \'etale
cover of $\cM_g$, whose degree is the
cardinality of $\PP^{2g}(\ZZ/n\ZZ)$. The space
$\cR^{(2)}\cR^{(2)}_g$ of unramified
$2$\/-sheeted covers of unramified $2$\/-sheeted
covers of curves of genus $g$ is therefore
\'etale over $\cR_g^{(2)}$. It contains two
subcovers of $\cR_g^{(2)}$, namely $\cR_g^{(4)}$
and the subspace $\cR_\sigma$ of
$\cR^{(2)}\cR^{(2)}_g$ determined by the
condition that $\Nm(\sigma^\ast) = \sigma$. The
Lemma asserts that these subcovers are equal,
and the ``only if'' direction just proved shows
the inclusion of $\cR_g^{(4)}$ into
$\cR_\sigma$. To prove equality it suffices to
show equality of the degrees of these two
subcovers over $\cR_g^{(2)}$. The degree of
$\cR_g^{(4)} /\cR_g^{(2)}$ is
$|\PP^{2g}(\ZZ/4\ZZ)/\PP^{2g}(\ZZ/2\ZZ)| =
2^{2g-1}$. The degree of $\cR_\sigma/
\cR_g^{(2)}$ is the cardinality of the kernel of
the norm map $\Nm$. This cardinality is the
quotient of the order of the source
$|H_1(F^\ast,\ZZ/2\ZZ)| = 2^{2(2g-1)}$ of $\Nm$
by the order of its image. As explained above,
this image is dual to the pull-back $\pi^\ast$
on homology. Hence this image has order
$|H_1(F,\ZZ/2\ZZ)|/|\Ker(\pi^\ast)| = 2^{2g}/2$,
so that the degree of $\cR_\sigma/ \cR_g^{(2)}$
is $2^{2g-1}$, which is the same as the degree
of $\cR_g^{(4)} / \cR_g^{(2)}$, concluding the
proof of the Lemma.
\end{proof}
\begin{proof}(of Lemma~\ref{l1})
Observe first that the even sets of Weierstrass
points of a hyperelliptic curve $V$ modulo
complementation are in natural bijection with
the points $\Jac(V)[2]$ of order $2$ of its
Jacobian. by Lemma~\ref{cyclic} the cover
$H^{\ast\ast}/H$ is cyclic if (and only if) the
norm map
\[\Nm:\Jac(H^\ast) \ra \Jac(H)\]
maps $\sigma^\ast$ to $\sigma$, where
$\sigma^\ast$ and $\sigma$ are the points of
order $2$ corresponding to $H^{\ast\ast}/H^\ast$
and $H^\ast/H$ respectively. Let
$q_1,\dots,q_{2g} \in \PP^1$ be the branch
points of $H/\PP^1$ other than $\infty,0$. Let
$H'$ be the double cover of $\PP^1$ branched
over $q_1,\dots,q_{2g}$, and let $\PP'$ be the
$u$\/-line, namely the double cover of $\PP^1$
branched above $\infty,0$. Let $q_i^+,q_i^-$ be
the inverse images of each $q_i$ in $\PP'$.
These are the branch points of the hyperelliptic
cover $H^\ast/\PP'$. Let $T$ be the
$\ZZ/2\ZZ$\/-module freely generated by the
$q_i$\/'s, Let $T'$ be the $\ZZ/2\ZZ$\/-module
freely generated by the $q_i^\pm$\/'s, and let
$t\in T$ and $t'\in T'$ be the sums of the
respective free generators. Let $T'_0$ be the
kernel of the natural degree map $\deg:T'\ra
\ZZ/2\ZZ$, and let $T^\ast$ be the quotient
$T'_0/\left<t'\right>$. Then we have natural
isomorphisms
\[ \Jac(H)[2] \simeq T \qquad\text{and}\qquad
\Jac(H^\ast)\simeq T^\ast.\] Under these
identifications $\sigma$ goes to $t$, while
$\Nm(q_i^+)=\Nm(q_i^-) = q_i$ for each $i$. It
follows that the covers $H^{\ast\ast}$\/'s which
are cyclic over $H$ correspond to those
$\sigma^\ast$\/'s which are a section of the
natural projection $\{q_i^+,q_i^-\}\ra \{q_i\}$,
namely to a choice of either $q_i^+$ or $q_i^-$
for each $i$. Let $H_1$ be the double cover of
$\PP'$ branched above the points of
$\sigma^\ast$, and let $H_2$ be the double cover
of $\PP'$ branched along the complementary
section. The following diagram shows the various
curves:
\[\begin{array}{rcccccl}
&&&&H^{\ast\ast}&&\\
&&&\swarrow&\downarrow\,&\searrow&\\
&&H^\ast&&H_1\,&&H_2\\
&\swarrow&\downarrow\,&\searrow&\downarrow\,&
\swarrow& \\
H'&&H&&\PP'&&\\
&\searrow&\downarrow&\swarrow&&&\\
&&\PP^1&&&&
\end{array}\]

It follows that $H^{\ast\ast}$ is the fibred
product $H_1\times_{\PP'} H_2$. At the same time
it is clear that the involution $u\mapsto -u$ of
$\PP'/\PP^1$ exchanges the branch loci of $H_1$
and $H_2$, so that these curves are isomorphic
to each other (and are otherwise arbitrary for
their fixed genus). Therefore
\[\Prym(H^{\ast\ast}/H^\ast) \simeq
\Jac(H_1)\times\Jac(H_2) \simeq \Jac(H_1)^2,\]
concluding the proof of the Lemma.

\begin{remark}
One can show that $\Prym(H^{\ast\ast}/H^\ast)$
is in fact {\em isomorphic} to
$\Jac(H_1)\times\Jac(H_2)$ using the bigonal
construction (see~\cite{don}). We omit the proof
since we do not need this.
\end{remark}
\end{proof}
\begin{proof} (of Lemma~\ref{l2})
The statement is well known when the genus
$g(C)$ of $C$ is $1$, so we assume by induction
that we know the result up to genus $g-1$ and
prove it for $g(C)=g>1$. We have a canonical
embedding
\[\End\Jac C\hookrightarrow \End H^1(\Jac C,\ZZ)
 = \End H^1(C,\ZZ). \]

Consider a degeneration of $C$ to a one point
union $C_0=C_1\cup C_2$, where each $C_i$ is a
general hyperelliptic curve of genus $g_i>0$ (if
$g_1=g_2$ we also assume that the Jacobians of
$C_1$ and of $C_2$ are not isogenous). Then
$\Jac C_0= \Jac C_1 \times \Jac C_2$. The family
of cohomology groups $H^1(C,\ZZ)$ and the Hodge
structures on $H^1(C)$ are continuous at $C_0$.
Hence $\End \Jac C$ embeds into $\End \Jac C_0$.
By our assumptions and by induction, we get a
decomposition $\End \Jac C_0\simeq \ZZ\oplus
\ZZ$ corresponding to the decomposition
$H^1(C_0,\ZZ)= H^1(C_1,\ZZ)\oplus H^1(C_2,\ZZ)$.
Since it is possible to make the degeneration of
$C$ to $C_0$ so as to give different such
decompositions, and since $\End \Jac C$ must
respect them all, it follows that $\End \Jac
C=\ZZ$ as asserted.
\end{proof}

\end{document}